\documentclass{etds}

\usepackage{amssymb}
\usepackage{color}

\newcommand{\area}{\operatorname{area}}

\newcommand{\cX}{\mathcal{X}}
\newcommand{\cM}{\mathcal{M}}
\newcommand{\cH}{\mathcal{H}}
\newcommand{\cO}{\mathcal{O}}
\newcommand{\cP}{\mathcal{P}}
\newcommand{\sing}{\Sigma}
\newcommand{\Xh}{\cX^k_{\hor}}
\newcommand{\Xv}{\cX^k_{\ver}}
\newcommand{\withsing}{\widehat}
\newcommand{\hol}{\operatorname{hol}}
\newcommand{\Vol}{\operatorname{Vol}}
\newcommand{\half}{\mathbb{H}}

\newcommand{\psgrp}{\Phi}
\newcommand{\semigrp}{\withsing{\Phi}}

\newcommand{\real}{\mathbb{R}}
\newcommand{\integer}{\mathbb{Z}}
\renewcommand{\natural}{\mathbb{N}}

\newcommand{\ga}{\mathsf{a}}
\newcommand{\gb}{\mathsf{b}}
\newcommand{\gc}{\mathsf{c}}
\newcommand{\gd}{\mathsf{d}}

\DeclareMathOperator{\Ad}{Ad}
\DeclareMathOperator{\SL}{SL}
\DeclareMathOperator{\SO}{SO}
\DeclareMathOperator{\Stab}{Stab}
\DeclareMathOperator{\VF}{Vect}
\DeclareMathOperator{\hor}{Horiz}
\DeclareMathOperator{\ver}{Vert}
\DeclareMathOperator{\Id}{Id}

\DeclareMathOperator{\domain}{domain}

\newcommand{\Lie}[1]{\mathfrak{\lowercase{#1}}}

\newcommand{\bigset}[2]{\left\{\, #1
      \mathrel{\left| \vphantom {\left\{ #1 \mid #2 \right\} }
      \right.} #2 \,\right\} }

\newcommand{\pref}[1]{{\upshape(\ref{#1})}}
\newcommand{\see}[1]{{\upshape(}see~\ref{#1}{\upshape)}}
\newcommand{\seeSect}[1]{{\upshape(}see~\S\ref{#1}{\upshape)}}
\newcommand{\fullref}[2]{\ref{#1}\pref{#1-#2}}
\newcommand{\fullsee}[2]{(see~\ref{#1}\pref{#1-#2})}

\numberwithin{equation}{section}

\newtheorem{thm}[equation]{Theorem}
\newtheorem{prop}[equation]{Proposition}
\newtheorem{lem}[equation]{Lemma}
\newtheorem{cor}[equation]{Corollary}

\newenvironment{rmk}{\refstepcounter{equation}\proc{Remark
\theequation.}}{\medbreak}

\newenvironment{notation}{\refstepcounter{equation}\proc{Notation
\theequation.}}{\medbreak}

\newenvironment{eg}{\refstepcounter{equation}\proc{Example
\theequation.}}{\medbreak}

\newenvironment{assump}{\refstepcounter{equation}\proc{Assumption
\theequation.}}{\medbreak}

\newenvironment{proof}{\proc{Proof.}}{\ep\medbreak}
\newenvironment{proofsketch}{\proc{Proof {\rm
(sketch)}.}}{\ep\medbreak}

\newenvironment{thmref}{\thmrefer}{}
\newcommand{\thmrefer}[1]{\renewcommand\theequation
       {\protect\ref{#1}$'$}\addtocounter{equation}{-1}}

\begin{document}

\ETDS{1}{36}{???}{2005} 

  \runningheads{A.~Eskin, J.~Marklof, and D.~Morris}{Unipotent flows and
Veech surfaces}

  \title{Unipotent flows on the space of branched covers of Veech 
surfaces}

  \author{ALEX ESKIN\affil{1}, JENS MARKLOF\affil{2}\ and
  DAVE WITTE MORRIS\affil{3}}
  \address{\affilnum{1}\
  Department of Mathematics, University of Chicago,
  Chicago, IL 60637 USA \\
  \email{eskin@math.uchicago.edu}\\
  \affilnum{2}\
  School of Mathematics, University of Bristol,
  Bristol BS8~1TW, U.K.\\
  \email{j.marklof@bristol.ac.uk}\\
  \affilnum{3}\
  Department of Mathematics and Computer Science, University of
Lethbridge,
  Lethbridge, Alberta, T1K~3M4, Canada \\
  \email{Dave.Morris@uleth.ca}}

\recd{July 2004; Revised May 2005}  

  \begin{abstract}
  There is a natural action of $\SL(2,\real)$ on the moduli space of
translation surfaces, and this yields an action of the unipotent
subgroup
  $U = \left\{ \begin{pmatrix} 1 & * \\ 0 & 1 \end{pmatrix} \right\}$.
  We classify the $U$-invariant ergodic measures on certain special
submanifolds of the moduli space. (Each submanifold is the
$\SL(2,\real)$-orbit of the set of branched covers of a fixed Veech
surface.)
  For the $U$-action on these submanifolds, this is an analogue of
Ratner's Theorem on unipotent flows.
  The result yields an asymptotic estimate of the number of periodic
trajectories for billiards in a certain family of non-Veech rational
triangles, namely, the isosceles triangles in which exactly one angle
is $2 \pi/n$, with $n \ge 5$ and $n$~odd.
  \end{abstract}

\section{Introduction}

  A polygon $P \subset \real^2$ is called rational if all angles of $P$
are rational multiples of $\pi$. Let $N(P,T)$ denote the number of
(cylinders of) periodic billiard trajectories of Euclidean length at
most $T$. It is a theorem of H.~Masur \cite{Masur:upper, Masur:lower}
that there exist constants $c_1 = c_1(P)$ and $c_2 = c_2(P)>0 $ such
that for $T \gg 1$,
  \begin{equation}
  \label{eq:masur}
  c_1 T^2 < N(P,T) < c_2 T^2.
  \end{equation}
  A natural question is whether (\ref{eq:masur}) can be converted to an
asymptotic formula as $T \to \infty$.

A well known construction associates a ``translation surface'' $S$ to
each rational polygon $P$. Essentially the algorithm ``unfolds'' the
billiard trajectories, by reflecting the polygon instead of reflecting
the trajectory. More precisely, let $\Delta \subset O(2)$ denote the
group generated by reflections in the sides of the polygon $P$. Since
$P$ is rational, $\Delta$ is finite. The ``translation surface''
consists of $\Delta$ copies of $P$, with each copy glued to each of
its mirror images along the reflecting side.

For example, if $P$ is the unit square, then $S$ is the torus
$\real^2/2 \integer \oplus 2 \integer$, and if $P$ is the isosceles
triangle with angles $\pi/2 - \pi/n$, $\pi/2 - \pi/n$, $2\pi/n$, and
$n$ is even, then $S$ is the regular $n$-gon with opposite sides
identified.

A translation surface can be defined in one of the following
equivalent ways:
\begin{enumerate} \renewcommand{\theenumi}{\alph{enumi}}
\item \label{TransSurfDefn-polygons}
  A union of polygons $P_1 \cup \dots \cup P_n$ where
      each $P_i \subset \real^2$, and the $P_i$ are glued along parallel
      sides, such that each side is glued to exactly one other, and the
      total angle in each vertex is an integer multiple of $2 \pi$.
\item An orientable surface with a flat metric and isolated
      conical singularities that has trivial rotational holonomy. (Note
that
      trivial rotational holonomy means in particular that parallel
transport
      of a vector along a small loop going around a conical point brings 
a
      vector back to itself. This implies that all cone angles are 
integer
      multiples of $2 \pi$.)
\item A pair $(M,\omega)$, where $M$ is an (orientable)
      Riemann surface, and $\omega$ is a holomorphic 1-form on $M$.
     (Note that away from the zeroes of $\omega$, there is a local
     coordinate $z$ such that $\omega = dz$, and this coordinate is 
unique
     up to translation. Then one can define the metric on $M$ as 
$|dz|^2$.
     This metric is flat, with conical singularities appearing at the
     zeroes of $\omega$.
\end{enumerate}
The term ``translation surface'' comes from the fact that away from
the cone points the surface can be covered by charts so that the
transition functions are translations ($z \to z+c$). If $\alpha =
(\alpha_1, \dots, \alpha_k)$ is an $n$-tuple of positive integers
such that the sum of the $\alpha_i$ is even, we denote by
$\cH(\alpha)$ the moduli space of translation surfaces $(M,\omega)$
such that the multiplicities of the zeroes of $\omega$ are given by
$\alpha_1, \dots, \alpha_n$ (or equivalently such that the orders of
the conical singularities are $2 \pi (\alpha_1 + 1), \dots, 2
\pi(\alpha_n + 1)$). (Actually, for technical reasons, the
singularities of $(M,\omega)$ should be labeled; thus, an element of
$\cH(\alpha)$ is a tuple $(M,\omega, p_1,\ldots,p_n)$, where
$p_1,\ldots,p_n$ are the singularities of~$M$, and the multiplicity
of~$p_i$ is~$\alpha_i$.) The moduli space of translation surfaces is
naturally stratified by the spaces $\cH(\alpha)$; each is called a
{\em stratum}.

By construction, billiard trajectories on $P$ correspond to ``straight
lines'' on $S$, which are geodesics not passing through singularities.
It is easy to see that any such geodesic
is part of a family of freely homotopic parallel geodesics of the
same length. Such a family is called a cylinder. Let
$N(S,T)$ denote the number of cylinders on $S$ of length at most $T$.
(By the length of a cylinder we mean the length of any of the closed
geodesics that comprise it).

\subsection*{The $\SL(2,\real)$ action.}
There is an action of $\SL(2,\real)$ on the moduli space of
translation surfaces that preserves the stratification. For our
purpose, it is easiest to see this using definition
\pref{TransSurfDefn-polygons}: since $\SL(2,\real)$ acts on $\real^2$,
for $S = P_1 \cup \dots \cup P_n$, we can define $g S = g P_1 \cup \dots
g P_n$, where all identifications between the sides of the polygons for
$gS$ are the same as for $S$. This action generalizes the action of
$\SL(2,\real)$ on the space of flat tori $\SL(2,\real)/\SL(2,\integer)$.

We can visualize this as a composition of ``the
usual linear action'' with ``cut and paste.'' We note that ``cut and
paste'' is an isometry on the surface (and in fact preserves the
horizontal and vertical directions as well).
Note that if $S$ is a union of triangles, and $g$ is a large element
of $\SL(2,\real)$ then $gS$ is a union of long and thin triangles.
We may if we wish ``cut and paste'' $gS$ and retriangulate to try to
present $gS$ as a union of triangles with bounded side lengths.

\subsection*{Veech surfaces.}
For $S \in \cH(\alpha)$, let $\Gamma(S) \subset \SL(2,\real)$
denote the stabilizer of $S$. The group $\Gamma(S)$ is called the
{\em Veech group} of $S$. If $\Gamma(S)$ is a lattice in $\SL(2,\real)$
then $S$ is called a {\em Veech surface}. It is a theorem of
Veech \cite{Veech:eisenstein} that if $S$ is a Veech surface, then
there
exists $c = c(S)$ such that
\begin{equation}
\label{eq:asymp}
N(S,T) \sim c T^2
\end{equation}
as $T \to \infty$.

\subsection*{Counting and Ratner's Theorem.}
One has the formula \cite{Veech:Siegel}, (reproduced in \cite{EM})
\begin{equation}
\label{eq:count}
N(S,T) - N(S,T/2) \approx T^2 \int_0^{2 \pi} \hat{f}( a_t r_\theta S)
\, d\theta,
\end{equation}
where $a_t = \begin{pmatrix} e^t & 0 \\ 0 & e^{-t} \end{pmatrix}$,
$r_\theta = \begin{pmatrix} \cos\theta & \sin \theta \\ -\sin \theta &
      \cos \theta \end{pmatrix}$, and $t = \log T$. The left hand side
counts (cylinders) of closed geodesics in an annulus, and the right
hand side is an integral over part of the $\SL(2,\real)$ orbit
of $S$. Thus, the $\SL(2,\real)$ action can be used to count closed
geodesics (and thus periodic billiard trajectories).

A closer examination of (\ref{eq:count}) shows that the integral is
over large circles inside the $\SL(2,\real)$ orbit. These large circles
can be approximated by horocycles, which are orbits of $u_t =
\begin{pmatrix} 1 & t \\ 0 & 1 \end{pmatrix}$. Thus the ergodic
properties of the action of $U = \{ u_t \; | \; t \in \real \}$
play a key role.

Ratner's theorem \cite{RatnerMeas}
is the classification of the invariant measures
for the action of a unipotent subgroup
on the homogeneous space $H/\Gamma$, where $H$ is a Lie group and
$\Gamma$ is a lattice in $H$.
An important question is whether a similar theorem holds for the
$U$-action on a stratum $\cH(\alpha)$. One can also ask this question
when one restricts the action to any $\SL(2,\real)$ invariant
submanifold of a stratum. In this paper, we will classify the
$U$-invariant measures on a certain family of $\SL(2,\real)$-invariant
manifolds. Another result in this direction was obtained by McMullen
\cite{McMullen:SL2}
who, in genus $2$, classified the measures invariant under all of
$\SL(2,\real)$.

\subsection*{Branched covers of Veech surfaces.}
We say that a translation surface $S$ is a branched cover of
a translation surface $M$ if the covering map $\pi$ respects the
translation structure (i.e. if we identify $S = (L_1, \omega_1)$ and
$M = (L_2, \omega_2)$ where the $L_i$ are Riemann surfaces and
the $\omega_i$ are holomorphic 1-forms on $L_i$ then we require
that $\pi\colon L_1 \to L_2$ is holomorphic and $\pi^*(\omega_2) =
\omega_1$.)

Now let $M \in \cH(\alpha)$ be
a Veech surface. Then the $\SL(2,\real)$ orbit of $M$ is a closed
subset $D$ of $\cH(\alpha)$.
Let $\cH(\beta)$ be another stratum, and
let
$\cM_D(\beta)$ denote the set of all translations surfaces $S \in
\cH(\beta)$ that are branched covers of $M \in D$.
We will always assume that $\beta$ is such that $\cM_D(\beta)$ is
not-empty. Then $\cM_D(\beta)$ is $\SL(2,\real)$ invariant.

There are two types of Veech surfaces: arithmetic and non-arithmetic.
A surface $S = (M,\omega)$ is an arithmetic Veech surface if and
only if $M$ is a (holomorphic) branched cover of a torus,  $\omega$
is the pullback by the covering map of the standard differential
$dz$ on the torus, and the branch points project to points of finite
order (under the additive group of the torus). Equivalently (see
\cite{Gutkin:Judge}) $S$ is an arithmetic Veech surface if and only
if $\Gamma(S)$ is commensurable to $\SL(2,\integer)$. All other
Veech surfaces are called non-arithmetic (and their Veech groups,
which are always non-uniform lattices, are non-arithmetic lattices
in $\SL(2,\real)$).
  The case where $M$ is arithmetic was analyzed in \cite{EMS}.

In this paper, we assume that $M$ is not arithmetic, which implies that
the genus of $M$ is greater then 1.
Then considering the Euler characteristic, it is easy to see that
the degree of $\pi$ is determined by $D$ and $\beta$. This implies
that $\cM_D(\beta)$ is closed. (In the case where the genus of $M$ is 1,
one also has to fix the degree of the cover; see \cite{EMS}
for the details.)

The main result of this paper is a classification of the $U$-invariant
ergodic measures on $\cM_D(\beta)$. This allows us to prove asymptotic
formulas of the form (\ref{eq:asymp}) for $S \in \cM_D(\beta)$ (see
Theorem~\ref{theorem:branched:cover:count}). In
particular we prove the following:

  \begin{thm}
  \label{theorem:triangle:count}
  Let $P_n$ be a triangle with angles
  $$ \frac{n-2}{2n} \pi, \ \frac{n-2}{2n} \pi, \ \frac{4}{2n} \pi
  ,$$ where $n \ge 5$, $n$ odd. Then, as $T \to \infty$,
  \begin{displaymath}
  N(P_n,T) \sim \frac{\pi}{\zeta(2)} \frac{(n-1)(n^2+n+3)}{144(n-2)}
\frac{T^2}{\area(P_n)} .
  \end{displaymath}
  \end{thm}

The fact that the surface $S_n$ associated to $P_n$ is not Veech but
is a branched cover of degree 2 of a Veech surface is due to P.~Hubert
and T.~Schmidt (see Proposition 4 in \cite{HS1} and its proof).
We should also note that if $n=5$ then the Veech group of $S_n$ is
infinitely generated (see \cite{HS2}). However, the Veech group of $S_n$
plays no direct role in our analysis.

Here is an outline of the paper.
  Section~\ref{MeasClassSect} states our main theorem.
  Section~\ref{PreliminariesSect} establishes notation and presents a few
basic lemmas.
  Section~\ref{ShearingSect} explains ``shearing," the foundation of our
study of invariant measures.
  Section~\ref{MainPfSect} proves our main theorem
\pref{theorem:classification} that classifies $U$-invariant measures.
  Section~\ref{CountabilitySection} proves that there are only countably
many closed orbits of a certain type.
  Section~\ref{AvgCircleSect} uses our main theorem (and the countability
result of \S\ref{CountabilitySection}) to prove that large circles in
$\SL(2,\real)$-orbits become uniformly distributed with respect to
certain natural measures.
  Section~\ref{Application} applies the equidistribution result of
\S\ref{AvgCircleSect} to derive asymptotic estimates for the number of
periodic trajectories in branched covers of Veech surfaces.

\section{Measure classification}
   \label{MeasClassSect}

\subsection*{Definitions and notation.}
   Let $G= \SL(2,\real)$.
   Let $M$ be a Veech surface, which means that $\Gamma = \Stab_G
\bigl( M \bigr)$ is a lattice in~$G$. Here, we use $M$ to also denote
the isometry class of~$M$; this is a single point in the moduli space.
For $k \in \natural$, we define $\cX^k$ to be the natural fiber bundle
over $G \cdot M$  whose fiber over $M$ is $M^k$. Thus, a point
of~$\cX^k$ is represented by $(M', p_1, \ldots, p_k)$, where $M' \in
GM$ and $p_1,\ldots,p_k \in M'$. In other words, a point in $\cX^k$
represents a surface in $M' \in G M$ together with $k$ marked points
on $M'$.

We note that the space $\cM_D(\beta)$ parameterizing branched covers
is itself a finite branched cover of the space $\cX^k$ for a suitable
$k$. (The covering map just maps $S \in M_D(\beta)$ to the surface in
$D$ it covers, and notes the locations of the branch points.) Thus,
to classify the $U$-invariant measures on $\cM_D(\beta)$ it is
enough to classify $U$-invariant measures on $\cX^k$
\see{FinCovClassify}.

If $M$ is a torus, then $\cX^k$ can be identified with the homogeneous
space $\bigl( G \ltimes (\real^2)^k \bigr) / \bigl( \SL(2,\integer)
\ltimes (\integer^2)^k \bigr)$. In this situation, a special case
of Ratner's Theorem \cite{RatnerMeas} classifies all the ergodic
$U$-invariant probability measures on~$\cX^k$. We generalize this to
allow $M$ to be any Veech surface. The proof is based heavily on
ideas of Ratner \cite{RatnerRig, RatnerQuot, RatnerJoin, RatnerSolv,
RatnerSS, RatnerMeas} and Margulis-Tomanov
\cite{MargulisTomanov-Ratner}. An introduction to these ideas can be
found in \cite{Morris-Ratner}.

Let $\sing$ be the singular set of~$M$. Then for $g \in G$, $g \sing$
is the singular set of $g M$. Let $M_0 = M \smallsetminus \sing$, and
let $\cX^k_0 \subset \cX^k$ denote the set $(gM, p_1, \ldots, p_k)$
where $g \in G$ and $\{p_1, \dots, p_k\} \cap g \sing = \emptyset$.
Then $\cX_0^k$ is isomorphic to the natural fiber bundle over $G  M $
whose fiber over $M$ is $(M_0)^k$.

We have a natural embedding of~$\real^2$ in the space $\VF(M_0)$ of
smooth vector fields on~$M_0$, so, for each $v \in \real^2$ and $p \in
M_0$, we have a trajectory $\gamma_{v,p}(t)$ in~$M_0$ that is defined
for $t$ in a certain open interval containing~$0$ (until the
trajectory hits the singular set). We are interested only in the
forward trajectory, that is, for $t \ge 0$. By including the singular
points of~$M$, we extend $\gamma_{v,p}$ to a continuous curve
$\withsing{\gamma}_{v,p}$ in~$M$ that is defined for $t$ in a closed
interval (and for all points in~$M$):
    \begin{itemize}
   \item let $\withsing{\gamma}_{v,p}(0) = p$ for all $v \in \real^2$
and $p \in M$; and
    \item if $t > 0$ and $t$ is in the closure of the domain of
$\gamma_{v,p}$, let
    $$\withsing{\gamma}_{v,p}(t) = \lim_{s \to t^-} \gamma_{v,p}(s) \in
\sing .$$
    \end{itemize}
    Then each $v \in \real^2$ defines a function $\withsing{\phi}_v
\colon  M_v \to M$, defined by $\withsing{\phi}_v(p) =
\withsing{\gamma}_{p,v}(1)$, where $M_v$ is a dense, open subset
of~$M$. Note that $\withsing{\phi}_v$ is a local isometry (hence
continuous). On the other hand, $\withsing{\phi}_v$ is usually not
invertible, because a singular point will typically have several
preimages. In addition, $\withsing{\phi}_v$ is usually not uniformly
continuous, because of branch cuts.

\begin{figure}
  \begin{center}
  \input{opposite.pstex_t}
  \caption{In our notation, $v \in \real^2$ and $w \in \real^2$ can be
close, but $\withsing{\phi}_v(p)$ and $\withsing{\phi}_w(p)$ may not be
close. The wavy line represents a branch cut.} \label{fig:opposite}
  \end{center}
  \end{figure}

For $w \in (\real^2)^k$, we have a continuous map $\withsing{\phi}^k_w
\colon \cX^k_w \to \cX^k$ (where $\cX^k_w$ is a certain subset
of~$\cX^k$), defined by
    $$ \withsing{\phi}^k_w (M,p_1,\ldots,p_k) = \bigl(M,
\withsing{\phi}_{v_1}(p_1),\ldots, \withsing{\phi}_{v_k}(p_k) \bigr)
.$$
  (Thus, $\withsing{\phi}^k_w$ does not change the surface~$M$, but
moves the marked points in the directions specified by~$w$.)
    Let $\semigrp^k_{(\real^2)^k}$ be the pseudosemigroup generated by
$\{\, \withsing{\phi}^k_w \mid w \in (\real^2)^k\,\}$. (The prefix
``pseudo" simply refers to the fact that these maps are not defined
on the entire space $\cX^k$, but only on a subset.) Although the maps in
$\semigrp^k_{(\real^2)^k}$ may not be one-to-one, they are always
finite-to-one.

For $w \in (\real^2)^k$, let $\phi^k_w$ be the restriction
of~$\withsing{\phi}^k_w$ to $(\withsing{\phi}^k_w)^{-1}(\cX^k_0)$.
Then $\phi^k_w$ is a diffeomorphism (and local isometry) from a dense
open subset of~$\cX_0$ to a dense open subset of~$\cX_0$. Let
$\psgrp^k_{(\real^2)^k}$ be the pseudogroup that is generated by $\{\,
\phi^k_w \mid w \in (\real^2)^k \,\}$. We remark that
$\psgrp^k_{(\real^2)^k}$ is transitive on~$\cX^k_0$.

    Note that each of $\psgrp_{(\real^2)^k}^k$ and
$\semigrp_{(\real^2)^k}^k$ is normalized by the action of $G =
\SL(2,\real)$ on~$\cX^k$, so we have corresponding semidirect
products $G \ltimes \psgrp_{(\real^2)^k}^k$ and $G \ltimes
\semigrp_{(\real^2)^k}^k$.

Let
    $$\hor = \{\, \bigl( (x_i,0) \bigr)_{i=1}^k \mid x_i \in \real \,\}
\subset (\real^2)^k $$
    and
    $$\withsing{\hor} = \{\, \withsing{\phi}^k_w \mid w \in \hor \,\} .$$
    Note that, for $w_1,w_2 \in \hor$, we have
$\withsing{\phi}^k_{w_1+w_2} = \withsing{\phi}^k_{w_1}
\withsing{\phi}^k_{w_2}$ on the intersection of their domains, so
$\withsing{\hor}$ is a pseudosemigroup. Also, $\withsing{\hor}$ commutes
with the action of $U$.

\subsection*{Statement of the main results.} Let $\mu$ be an ergodic
$U$-invariant probability measure on~$\cX^k$. The projection of~$\mu$
to $G/\Gamma$ is $U$-invariant, so it must be either Lebesgue measure
or the arc-length on a closed $U$-orbit \cite{Dani-RatnerSL2}. The
interesting case is when the projection is Lebesgue. A weak statement
of our results is simply to say that, in this case, some horizontal
translate of~$\mu$ must be $G$-invariant:

\begin{thm} \label{theorem:U->G}
      Suppose $\mu$ is any ergodic $U$-invariant probability measure on
$\cX^k$, such that the projection of~$\mu$ to $G/\Gamma$ is Lebesgue.
Then there exists $h \in \withsing{\hor}$, such that $h_*\mu$ is
$G$-invariant {\rm(}and the domain of~$h$ has full measure{\rm)}.
      \end{thm}

To obtain a more precise description of the $U$-invariant measures,
one need only describe the $G$-invariant measures on $\cX^k$.

\begin{rmk} \label{Remark:BeforeThm} \
      \begin{enumerate}
      \item \label{Remark:BeforeThm-Gamma=G}
      It is easy to see that the $G$-invariant probability measures
on~$\cX^k$ are in natural one-to-one correspondence with the
$\Gamma$-invariant probability measures on~$M^k$ (cf., e.g.,
\cite[pf.\ of Cor.~5.8]{Witte:Quotients}).
      \item \label{Remark:BeforeThm-NoSing}
    It is the $\Gamma$-invariant measures on $M^k_0$ that are the most
important to understand, because it is easy to see that every ergodic
measure on $M^k$ arises from the following construction.
    Choose some $p_1 \in \sing^d$ and some probability measure~$\nu$ on
$M^{k-d}_0$ that is invariant under a finite-index subgroup
of~$\Gamma$. The corresponding measure on $\{p_1\} \times M^{k-d}$ is
invariant under a finite-index subgroup~$\Gamma'$ of~$\Gamma$. By
averaging over $\Gamma/\Gamma'$, this yields a $\Gamma$-invariant
measure supported on the subset $\sing^d \times M^{k-d}$ of $M^k$.
      \end{enumerate}
      \end{rmk}

We will show that every ergodic measure is carried by a nice subspace
of~$M^k$. In particular, any ergodic measure carried by $M^k_0$ is the
Lebesgue measure on a flat submanifold of $M^k_0$.

\begin{eg}
      The natural Lebesgue measure on the diagonal $\Delta =
\{(p,p,p)\}$ of $M^3_0$ is a $\Gamma$-invariant probability measure
on $M^3_0$. Note that $W = \{v,v,v\}$ is a $G$-invariant subspace of
$(\real^2)^3$, and that the pseudogroup $\psgrp_W^3$ of
diffeomorphisms it generates is transitive on~$\Delta$.
      \end{eg}

\begin{thm} \label{theorem:GammaInvt}
      Suppose $\mu$ is an ergodic $\Gamma$-invariant probability
measure on $M^k$. Then there exist
      \begin{itemize}
      \item a point $p \in M^k$,
      and
      \item a $G$-invariant linear subspace~$W$ of $(\real^2)^k$,
      \end{itemize}
      such that
      \begin{enumerate}
      \item the orbit $\semigrp_W^k(p)$ of~$p$ under $\semigrp_W^k$ is a
closed subset of~$M^k$ whose dimension is $\dim W$,
      \item \label{theorem:GammaInvt-finiteindex}
  some finite-index subgroup of\/~$\Gamma$ fixes $\semigrp_W^k(p)$
setwise,
      and
      \item $\mu$ is the $\semigrp_W^k(p)$-invariant Lebesgue measure on
$\Gamma \semigrp_W^k(p)$.
      \end{enumerate}
      \end{thm}

\begin{rmk} \
  \begin{enumerate}
      \item Conversely, if $W$ is $G$-invariant, $\semigrp_W^k(p)$ is
closed, and some finite-index subgroup of~$\Gamma$ fixes
$\semigrp_W^k(p)$, then the $\semigrp_W^k$-invariant Lebesgue measure on
$\Gamma \semigrp_W^k(p)$ is a $\Gamma$-invariant probability measure.
However, it may not be ergodic.
      \item We wish to emphasize that
conclusion~\fullref{theorem:GammaInvt}{finiteindex} implies the set
$\Gamma \semigrp_W^k(p)$ is a finite union of translates of
$\semigrp_W^k(p)$.
      \end{enumerate}
      \end{rmk}

The theorem can be stated in the following equivalent form
\fullsee{Remark:BeforeThm}{Gamma=G}:

\begin{thmref}{theorem:GammaInvt} \begin{thm}
\label{theorem:GInvtMeas}
      Suppose $\mu$ is an ergodic $G$-invariant probability measure on
$\cX^k_0$. Then there exist
      \begin{itemize}
      \item a point $(M,p) \in \cX^k$
      and
      \item a $G$-invariant linear subspace~$W$ of $(\real^2)^k$,
      \end{itemize}
      such that
      \begin{enumerate}
      \item the orbit $(G \ltimes \semigrp_W^k) p$ of~$p$ under $G
\ltimes \semigrp_W^k$ is a closed subset of~$\cX^k$ whose dimension is
$\dim(G \ltimes \semigrp_W^k)$,
      and
      \item $\mu$ is the $(G \ltimes \semigrp_W^k)$-invariant Lebesgue
measure on this orbit.
      \end{enumerate}
      \end{thm}
      \end{thmref}

This results in the following explicit version of
Theorem~\ref{theorem:U->G}:

      \begin{thm} \label{theorem:classification}
      Suppose $\mu$ is an ergodic $U$-invariant probability measure on
$\cX^k_0$. Then there exist
      \begin{itemize}
      \item a point $(M,p) \in \cX^k$,
      \item a $G$-invariant subspace~$W$ of $(\real^2)^k$,
      and
      \item some $h \in \hor$,
      \end{itemize}
      such that
      \begin{enumerate}
      \item $\mu \bigl( \domain(\withsing{\phi}^k_h) \bigr) = 1$,
      \item the orbit $(G \ltimes \semigrp_W^k) p$ of~$p$ under $G
\ltimes \semigrp_W^k$ is a closed subset of~$\cX^k$ whose dimension is
$\dim G + \dim W$,
      and
      \item $(\withsing{\phi}^k_h)_*\mu$ is the $(G \ltimes
\semigrp_W^k)$-invariant Lebesgue measure on this orbit.
      \end{enumerate}
      \end{thm}

We will give an application to counting the number of periodic
trajectories on~$M$ \seeSect{Application}.

Theorems~\ref{theorem:U->G}, \ref{theorem:GammaInvt},
and~\ref{theorem:GInvtMeas} have been stated only for expository
purposes --- they are not a part of the logical development. We prove
only Theorem~\ref{theorem:classification}, and the interested reader can
easily derive the other theorems as corollaries.

Our results imply that the closure of every $\Gamma$-orbit in~$M^k$ is
of a nice geometric form. Since $M \smallsetminus M_0 = \sing$ is a
$\Gamma$-invariant finite set, it suffices to describe the orbits of
points in $M_0^k$:

\begin{thmref}{corollary:GOrbit}
  \begin{cor} \label{corollary:GammaOrbit}
      Suppose $p \in M_0^k$. Then there exists a $G$-invariant linear
subspace~$W$ of $(\real^2)^k$,
      such that
      \begin{enumerate}
      \item the orbit $\semigrp_W^k(p)$ of~$p$ under $\semigrp_W^k$ is a
closed subset of~$M^k$ {\rm(}and its dimension is $\dim W${\rm)},
      \item some finite-index subgroup of~$\Gamma$ fixes
$\semigrp_W^k(p)$ setwise,
      and
      \item $\Gamma \semigrp_W^k(p)$ is the closure of the $\Gamma$-orbit
of~$p$.
      \end{enumerate}
      \end{cor}
  \end{thmref}

\section{Preliminaries} \label{PreliminariesSect}

We collect all the notation in this section. Some of this repeats the
definitions given in the previous sections.

\begin{notation} \label{PrelimNotation} \
      \begin{itemize}

      \item Let $G = \SL(2,\real)$.
      \item There is a natural action of~$G$ on the moduli space of
translation surfaces. We can visualize this as a composition of ``the
usual linear action'' with ``cut and paste.'' We note that ``cut and
paste'' is an isometry on the surface (and in fact preserves the
horizontal and vertical directions as well).
      \item Let $M$ be a Veech surface, which means that $\Gamma =
\Stab_G \bigl( M \bigr)$ is a lattice in~$G$. Here, we use $M$ to
also denote the isometry class of~$M$; this is a single point in the
moduli space.
      \item Let $k \in \natural$.
      \item We define $\cX^k$ to be the natural fiber bundle over $G M$
whose fiber over $ M$ is $M^k$. Thus, a point of~$\cX^k$ is
represented by $(M', p_1, \ldots, p_k)$, where $M' \in GM$ and
$p_1,\ldots,p_k \in M'$.
      \item The metric on $\cX^k$ is defined by
      \begin{align*}
      & d_{\cX^k} \bigl([M', (p_i)_{i=1}^k], [(M'', (q_i)_{i=1}^k]
\bigr) \\
      & \hskip0.5in = \min_{\begin{matrix}
      g \in G, \\ g M' = M''
      \end{matrix}}
      \left( \|g - \Id\| + \sum_{i=1}^k d_{M''}(g p_i,q_i) \right) .
      \end{align*}
      \item Note that $\cX^k$ is $G$-equivariantly homeomorphic to $(G
\times M^k)/\Gamma$, where
      \begin{itemize}
      \item $\Gamma$ acts on $G$ by right multiplication,
      \item $\Gamma$ acts on~$M^k$ componentwise,
      and
      \item $G$ acts on $\cX^k$ via $g \bigl( h, (p_i)_{i=1}^k \bigr) =
\bigl( gh, g(p_i)_{i=1}^k \bigr)$.
      \end{itemize}
      \item Let $\sing$ be the singular set of~$M$.
      \item Let $M_0 = M \smallsetminus \sing$.
      \item Let $\cX^k_0 = (G \times M_0^k)/\Gamma \subset \cX^k$.
      \item Any $w \in (\real^2)^k$ naturally defines a vector field
on~$\cX_0$. By taking the time-one map of the corresponding flow
(where it is defined), we obtain a diffeomorphism $\phi^k_w$ between
two dense open subsets of~$\cX_0$. The collection $\{\, \phi^k_w \mid
w \in (\real^2)^k \,\}$ generates a transitive pseudogroup
$\psgrp^k_{(\real^2)^k}$ of local diffeomorphisms of~$\cX^k_0$.

We extend $\phi^k_w$ to a (continuous) transformation
$\withsing{\phi}^k_w$ that is defined on a slightly larger subset
of~$\cX$, by letting
   $$ \withsing{\phi}^k_w(x)
   = \lim_{\begin{matrix} x' \to x \\ x' \in \domain \phi^k_w
\end{matrix}}
   \phi^k_w(x') $$
   if the limit exists. (See \S\ref{MeasClassSect} for a more concrete
definition of $\withsing{\phi}^k_w$, in terms of the flow
corresponding to~$w$.) We let $\semigrp_{(\real^2)^k}^k$ be the
pseudosemigroup generated by these maps.

Because the action of $G$ on $\cX^k$ normalizes
$\psgrp^k_{(\real^2)^k}$ and $\psgrp^k_{(\real^2)^k}$, we have
semidirect products $G \ltimes \psgrp_{(\real^2)^k}$ and $G \ltimes
\semigrp_{(\real^2)^k}$. Note that $G \ltimes \psgrp_{(\real^2)^k}$
is transitive on $\cX^k_0$.

It is important to note that, because of the singularities and
resulting branch cuts, $\phi^k_w$ is usually not uniformly continuous
(even though it is a local isometry). Furthermore, $\phi^k_w(p)$ is
not a uniformly continuous function of~$w$. See
Figure~\ref{fig:opposite}. Abusing notation, we may sometimes write
$w + p$ instead of $\phi_w(p)$.
      \item Let $U = \{\, u^t \mid t \in \real \,\}$, where $u^t =
\begin{bmatrix} 1 & t \\ 0 & 1 \end{bmatrix} \in G$.
      \item Let $A = \{\, a^s \mid s \in \real \,\}$, where $a^s =
\begin{bmatrix} e^s & 0 \\ 0 & e^{-s} \end{bmatrix} \in G$.
      \item Let $V = \{\, v^r \mid r \in \real \,\}$, where $v^r =
\begin{bmatrix} 1 & 0 \\ r & 1 \end{bmatrix} \in G$.
      \item Let $\mu$ be a $U$-invariant probability measure on $\cX^k$,
such that $\mu$ projects to the Lebesgue measure on $\Gamma
\backslash G$.
      \item Let $\hor = \bigset{\bigl( (x_i,0) \bigr)_{i=1}^k }{ x_i \in
\real }  \subset (\real^2)^k$
    and
     $\withsing{\hor} = \{\, \withsing{\phi}^k_w \mid w \in \hor \,\}$.
  Then $\withsing{\hor}$ is a pseudosemigroup.
      \item Let $\ver =
\{\, \bigl( (0,y_i) \bigr)_{i=1}^k \mid y_i \in \real \,\}  \subset
(\real^2)^k$ and
     $$\text{$\withsing{\ver}$ be the pseudosemigroup generated by $\{\,
\withsing{\phi}_w \mid w \in \ver \,\}$} .$$

      \item For $s \in \real$, we define $H_s \colon (\real^2)^k \to
\hor$ by
      $H_s(w) = u^s w - w$. Thus,
      $$ H_s \bigl( (x_i,y_i)_{i=1}^k \bigr)
      = (s y_i, 0)_{i=1}^k .$$

      \item The set
    $$\Xh = \{\, p \in \cX_0 \mid \withsing{\hor} \, p \subset \cX_0
\,\}  $$
    is $U$-invariant. Thus, it is either null or conull. Let us assume
it  is conull. (If not, then by ergodicity, there exists $h \in
\semigrp_{\hor}^k$ such that $h_*\mu$ is supported on $\cX
\smallsetminus \cX_0$. So $h_*\mu$ can be described by a construction
similar to Remark~\fullref{Remark:BeforeThm}{NoSing}. The conclusion of
Theorem~\ref{theorem:classification} is therefore obtained by induction
on~$k$.)

Note that $\hor$ acts on~$\Xh$, by $x(p) = \phi^k_x(p)$. Therefore,
the group $AU\ltimes \hor$ acts on $\Xh$.
    \item Let
    $$\Xv = \{\, p \in \cX_0 \mid \withsing{\ver} \, p \subset \cX_0
\,\}  $$
     Note that $AV\ltimes \ver$ acts on $\Xv$, but we do not yet know
that $\Xv$ is conull.
      \item Let $X = \{\, x \in \hor \mid x_* \mu = \mu \,\}$. Because
$\hor$ acts on $\Xh$, we know that $X$ is a closed subgroup of $\hor$.
      \item Let $Y = (v^1 - \Id)X \subset \ver$. Equivalently,
      $$ Y = \{\, y \in \ver \mid \mbox{$H_s(y) \in X$, for all $s \in
\real$} \,\} .$$
      \item Let $W = X + Y$. Note that $W$ is a $G$-invariant subspace
of~$(\real^2)^k$, so $G \ltimes \semigrp^k_W$ is a pseudosemigroup.
      \item Let $d = \dim X$.
      \item Let $\hor \ominus X = \hor \cap \bigl( 0^d \times
(\real^2)^{k-d} \bigr)$. By permuting coordinates, we may assume $X
\cap (\hor \ominus X) = 0$.
      \item Let $\pi_i \colon \cX^k \to \cX^i$ (the first $i$
coordinates) be the natural projection.
      \item For $\omega \in \cX^k$, we use $\mu_{\pi_i(\omega)}$ to
denote the fiber measure of~$\mu$ over the point $\pi_i(\omega)$
of~$\cX^i$.
      \end{itemize}
      \end{notation}

The following is obtained by applying the Pointwise Ergodic Theorem
to the action of~$U$ on~$\cX^k$.

\begin{lem}[cf.\ {\cite[Lem.~7.3]{MargulisTomanov-Ratner}}]
\label{UnifGeneric}
      For any $\rho > 0$, there is a ``uniformly generic set"
$\Omega_\rho$ in~$\cX^k$, such that
      \begin{enumerate}
      \item $\mu(\Omega_\rho) > 1-\rho$,
      \item for every $\epsilon > 0$ and every compact subset $K$
of~$\cX^k$, with $\mu(K) > 1 - \epsilon$, there exists $L_0 \in
\real^+$, such that, for all $\omega \in \Omega_\rho$ and all $L >
L_0$, we have
      $$ \lambda \{\, s \in [-L,L] \mid d(u^s \omega, K) < \epsilon
\,\} > (1 - \epsilon) (2L) ,$$
      where $\lambda$ is the Lebesgue measure on~$\real$.
      \end{enumerate}
      \end{lem}

\begin{lem}[cf.\ {\cite[Thm.~2.2]{RatnerSolv}},
{\cite[Lem.~5.8.6]{Morris-Ratner}}] \label{simultaneous}
      Suppose a Lie group~$H$ acts continuously on a Borel subset~$M$
of a locally compact metric space. If
      \begin{itemize}
      \item $U$ is a one-parameter, normal subgroup of~$H$,
      and
      \item $\mu$ is an ergodic $U$-invariant probability measure
on~$M$,
      \end{itemize}
      then
      \begin{enumerate}
      \item \label{simultaneous-conull}
   there is a $U$-invariant, Borel subset~$\Omega$ of~$M$, such that
      \begin{enumerate}
      \item $\mu(\Omega) = 1$,
      and
      \item $\Omega \cap c \Omega = \emptyset$ for all $c \in H
\smallsetminus \Stab_H(\mu)$,
      \end{enumerate}
      and
      \item \label{simultaneous-cpct}
   for any $\epsilon > 0$, there is a compact subset~$K$
of~$M$, such that
      \begin{enumerate}
      \item $\mu(K) > 1 - \epsilon$,
      and
      \item $K \cap c K = \emptyset$ for all $c \in H \smallsetminus
\Stab_H(\mu)$.
      \end{enumerate}
      \end{enumerate}
      \end{lem}

\begin{proof}
   Ratner's argument in \cite[Thm.~2.2]{RatnerSolv} shows, for each
$h_0 \in H \smallsetminus \Stab_H(\mu)$, that there is a
neighborhood~$B_{h_0}$ of~$h_0$ in $H \smallsetminus \Stab_H(\mu)$
and a conull $U$-invariant subset~$\Omega_{h_0}$ of~$M$, such that
      $$\mbox{$\Omega_{h_0} \cap h \Omega_{h_0} = \emptyset$, for all
$h  \in B_{h_0}$.}$$
   For the reader's convenience, we sketch the proof of this fact.
Because $h_0$ normalizes~$U$ but does not belong to $\Stab_H(\mu)$,
we know that $(h_0)_* \mu$ is $U$-invariant and ergodic, but is not
equal to~$\mu$. Therefore $(h_0)_* \mu$ and $\mu$ are mutually
singular, which implies there is a compact subset~$K_0$ of~$M$, such
that $\mu(K_0) > 0.99$ and $K_0 \cap h_0 K_0 = \emptyset$. By
continuity and compactness, there are open neighborhoods $\mathcal{U}$
and $\mathcal{U}^+$ of~$K_0$, and a symmetric neighborhood $B_e$
of~$e$ in~$H$, such that $\mathcal{U}^+ \cap h_0 (\mathcal{U}^+ \cap
M) = \emptyset$ and $B_e (\mathcal{U} \cap M) \subset \mathcal{U}^+$.
  From the Pointwise Ergodic Theorem, we know there is a conull
$U$-invariant subset~$\Omega_{h_0}$ of~$M$, such that the $U$-orbit
of every point in~$\Omega_{h_0}$ spends 99\% of its life in
$\mathcal{U} \cap M$. Now suppose there exists $h \in B_e h_0$, such
that $\Omega_{h_0} \cap h \Omega_{h_0} \neq \emptyset$. Then there
exists $x \in \Omega_{h_0}$, $u \in U$, and $c \in B_e$, such that $u
x$ and $c h_0 ux$ both belong to $\mathcal{U} \cap M$. This
implies that $ux$ and $h_0 ux$ both belong to~$\mathcal{U}^+$. This
contradicts the fact that $\mathcal{U}^+ \cap h_0 \mathcal{U}^+ =
\emptyset$.

      \pref{simultaneous-conull} Cover $H \smallsetminus
\Stab_H(\mu)$ with countably many balls $B_{h_j}$, and let $\Omega =
\bigcap_{j=1}^\infty \Omega_{h_j}$.

   \pref{simultaneous-cpct} Let $K$ be any compact subset
of~$\Omega$ with $\mu(K) > 1 - \epsilon$.
      \end{proof}

\begin{thm}[Kerckhoff-Masur-Smillie
{\cite[Thm.~2]{KerckhoffMasurSmillie}}] \label{KMS:UniqErg}
      For a.e.\ $v \in \real^2$, the foliation by orbits of $\real v$ is
uniquely ergodic on~$M_0$.
      \end{thm}

\begin{cor} \label{UniqErgHoriz}
      Suppose $\mu$ is a $U$-invariant probability measure on $\Xh$
whose projection to $G/\Gamma$ is Lebesgue.

      If $\mu$ is $\hor$-invariant, then $\mu$ is the Lebesgue measure.
      \end{cor}

\begin{proof}
      Theorem~\ref{KMS:UniqErg} implies that the foliation by orbits of
$\hor$ is uniquely ergodic on~$gM^k_0$, for a.e.\ $g \in G$. Thus,
almost every fiber of~$\mu$ over $G/\Gamma$ is the Lebesgue measure.
      \end{proof}

\section{Shearing} \label{ShearingSect}

In this section, we prove the crucial fact that the direction of fastest
transverse divergence between two nearby $U$-orbits is always along the
stabilizer of~$\mu$. The analogous statement for unipotent flows is a
cornerstone of the proof of Ratner's Theorem \cite[Lem.~3.3]{RatnerSS},
\cite[Lem.~7.5]{MargulisTomanov-Ratner},
\cite[Prop.~5.2.4$'$]{Morris-Ratner}.

\begin{notation} \label{DivNotn} \
      \begin{itemize}

      \item
      For any $g \in G$, we may write
      $$ g = \begin{bmatrix} 1+\ga & \gb \\ \gc & 1 + \gd \end{bmatrix} 
$$
      with $\ga,\gb,\gc,\gd \in \real$. For a sequence $\{g_n\} \subset
G$, we have $g_n \to e$ if and only if $\ga_n,\gb_n,\gc_n,\gd_n
\to 0$.

      \item Suppose $|\gd| < 1/4$, say. For $s \in \real$ with $|s| <
1/(4|\gc|)$, let
      \begin{itemize}
      \item $f(s,g) =  \displaystyle \frac{(1+\ga)s - \gb}{1 + \gd - \gc
s}
\in
\real$,
      \item $v_s(g) = \begin{bmatrix}
      1 & 0 \\
      (1 + \gd - \gc s) \gc & 1
      \end{bmatrix}
      \in V$,
      and
      \item $a_s(g) = \begin{bmatrix}
      1/(1 + \gd - \gc s) & 0 \\
      0 & 1 + \gd - \gc s
      \end{bmatrix}
      \in A$.
      \end{itemize}
     Note that $v_s(g) \to e$ if $g \to e$.

\item Suppose $\{p_n\}$ and $\{q_n\}$ are two sequences in a metric
space. If $d(p_n,q_n) \to 0$, we may write $p_n \approx q_n$.

      \end{itemize}
      \end{notation}

\begin{lem} \label{uguCalc}
        A simple calculation shows that
      $$
       u^{f(s,g)} g u^{-s} =
      \begin{bmatrix}
      1/(1 + \gd - \gc s) & 0 \\
      \gc & 1 + \gd - \gc s
      \end{bmatrix}
      = v_s(g) a_s(g).$$
     For a sequence $g_n \to e$, we denote $f_n(s_n) = f(s_n,g_n)$,
     and $a_{n,s_n} = a_{s_n}(g_n)$. Then
     $u^{f_n(s_n)} g_n u^{-s_n} \approx a_{n,s_n}$ if $g_n \to  e$
{\rm(}and $|s_n| < 1/(4|\gc_n|)${\rm)}.
     \end{lem}

\begin{rmk} (``Shearing")
      Let us discuss the action of~$U$ on~$(\real^2)^k$. For any $s \in
\real$ and $w \in (\real^2)^k$, we have
      $$u^s(w) = w + H_s(w) .$$
      Assume, now that
      $$ \mbox{$w_n,w'_n \to 0$ and $H_1(w_n) \neq H_1(w'_n)$.} $$
      There is some $s_n \in \real^+$, such that $\| H_{s_n}(w_n - w'_n)
\|
=
1$. Then
      $$ u^{s_n}(w_n) - u^{s_n}(w'_n) = (w_n - w'_n) + H_{s_n}(w_n - 
w'_n)
\approx H_{s_n}(w_n - w'_n) \in \hor .$$
      Thus, under the $U$-flow, $w_n$ and~$w'_n$ move apart along a leaf
of
the $\hor$-foliation. In other words, the direction in which two nearby
points move apart fastest is along~$\hor$.
      \end{rmk}

We use the notation of \pref{DivNotn} to state the main result of this
section.

      \begin{prop} \label{DivStab}
      For every $\rho > 0$, there is a compact subset $\Omega_\rho$
of~$\Xh$, with $\mu(\Omega_\rho) > 1 - \rho$, such that, if
      \begin{itemize}
      \item $(M_n,p_n), (M_n', p_n')$ are convergent sequences in
$\Omega_\rho$,
      \item $(M_n', p_n') = g_n w_n (M_n,p_n)$ for some $g_n \in G$ and
$w_n
\in  (\real^2)^k$,
      \item $g_n \to e$ and $w_n \to 0$,
     \item $s_n \in \real$ with
     $$|s_n| \le \frac{1}{\max \bigl( 4|\gc_n| , \|H_1(w_n)\| \bigr)} ,$$
     and
     \item $a_{n,s_n} H_{s_n}(w_n)$ converges,
      \end{itemize}
      then $\lim_{n \to \infty} a_{n,s_n} H_{s_n}(w_n) \in
\Stab_{A\hor}(\mu)^\circ$.
      \end{prop}

\begin{proof}
     Define
     $\varphi \colon [-1,1] \to A \hor$ by
     $$ \varphi(t)
     = \lim_{n \to \infty} a_{n,t|s_n|} H_{t|s_n|}(w_n)
     = a_\infty(t) \, h_\infty(t) ,$$
     where, letting $\gc_\infty = \lim_{n \to \infty} \gc_n |s_n|$ and
$w_\infty = \lim_{n \to \infty} |s_n| H_1(w_n)$, we have
     $$
     a_\infty(t) =
\begin{bmatrix}
      1/(1 - \gc_\infty t) & 0 \\
      0 & 1 - \gc_\infty t
      \end{bmatrix}
     \mbox{\qquad and\qquad}
     h_\infty(t) = t \, w_\infty .$$
     It is clear that $\varphi$ is continuous. We will show
$\varphi(t) \in \Stab_{A\hor}(\mu)$ for all~$t$. Then
     $$ \lim_{n \to \infty} a_{n,s_n} H_{s_n}(w_n) = \varphi(\pm 1)
\in \Stab_{A\hor}(\mu)^\circ ,$$
     as desired.

Let $\Omega_\rho$ be a uniformly generic set for the action of~$U$ on
$\cX^k$ with $\mu(\Omega_\rho) > 1 - \rho$ \see{UnifGeneric}. By
passing to a subset, we may assume that $\Omega_\rho \subset \Xh$ and
that $\Omega_\rho$ is compact. For any $\epsilon > 0$, we know, from
Lemma~\ref{simultaneous} (with $H = A U \hor$), that there is a compact
subset~$K$ of~$\Xh$, such that $\mu(K) > 1-(\epsilon/100)$ and $K
\cap h K = \emptyset$, for all $h \in A \hor \smallsetminus \Stab_{A
\hor}(\mu)$.

When $n$ is large, the definition of~$\Omega_\rho$ implies that
      \begin{equation} \label{DivStabPf-Uclose}
      d \bigl( u^s ( M_n,p) , K \bigr) < \epsilon
      \end{equation}
      for all but~$\epsilon$\% of the values of~$s$ in the interval
$\bigl[-|s_n|/4,|s_n|/4\bigr]$ (or longer intervals)
\see{UnifGeneric}. Note that the Jacobian of~$f_n$ is uniformly
bounded on $\bigl[ -|s_n|, |s_n| \bigr]$. More precisely, $f_n'(s) =
1/(1 + \gd_n - \gc_n s)^2$,  so $1/4 < f'(s) < 4$. Therefore,
      \begin{equation} \label{DivStabPf-U'close}
      d \bigl( u^{f(s)} (M_n',p') , K \bigr) < \epsilon
      \end{equation}
      for all but $4\epsilon$\% of the values of~$s$ in the interval
$\bigl[-|s_n|,|s_n|\bigr]$.
      Thus, \pref{DivStabPf-Uclose} and \pref{DivStabPf-U'close} hold
\emph{simultaneously} for all but $5\epsilon$\% of the values of~$s$ in
the interval $\bigl[-|s_n|,|s_n|\bigr]$.

      Let $(M,p) = \lim_{n \to \infty} (M_n,p_n)$. Because $(M,p) \in
\Omega_\rho \subset \Xh$, we know that translating $p_n$ by a
vector in~$\hor$ cannot move it into~$\sing$. Hence $d(C p, \sing) > 0$
for any compact subset~$C$ of $\hor$. Therefore, if $n$~is sufficiently
large, and, for convenience, we let $s = t |s_n|$, then
      \begin{itemize}
      \item $M_n$ has no singularities in
      $$ \bigset{ xyp_n }{ \begin{matrix}
      x \in \hor, y \in \ver, \\
       \|x\| \le 4\|H_{s_n}(w_n)\|, \\
      \|y\| \le 2 \|w_n\|
      \end{matrix} } ,$$
      so
      \item $u^sM_n$ has no singularities in
      $$  \bigset{ xy \, u^s p_n }{
      \begin{matrix}
      x \in \hor, y \in \ver, \\
       \|x\| \le 2\|H_{s_n}(w_n)\|, \\
      \|y\| \le 2 \|w_n\|
      \end{matrix} } .$$
      \end{itemize}
      This implies that
      $$ (u^s M_n, u^s p_n + u^s w_n) \approx \bigl( u^s M_n, u^s p_n +
H_s(w_n) \bigr)
.$$
      Therefore
      \begin{align*}
      u^{f_n(s)} (M_n',p_n')
      &= u^{f_n(s)} g_n u^{-s} (u^s M_n, u^s p_n + u^s w_n)
     && \text{(defn.\ of $g_n$ and $w_n$)} \\
      &\approx a_s \bigl( u^s M_n, u^s p_n + H_s(w_n) \bigr)
     && \text{(Lem.~\ref{uguCalc})} \\
      &= a_s H_s(w_n) u^s (M_n,p_n)
     && \text{($u^s$ commutes with $\hor$)} .
      \end{align*}
  When \pref{DivStabPf-Uclose} and \pref{DivStabPf-U'close} hold
simultaneously, we conclude that
  $$d \bigl( K, a_s H_s(w_n) K \bigr) \to 0 .$$
  From the definition of~$K$, we conclude that $a_s H_s(w_n) \in
\Stab_{A \hor}(\mu)$. That is, $\varphi(t) = a_s H_s(w_n)$ belongs to
$\Stab_{A \hor}(\mu)$ for all but $5\epsilon$\% of the values of~$t$ in
$[-1,1]$. Because $\epsilon$ is arbitrary, $\varphi$~is continuous, and
$\Stab_{A \hor}(\mu)$ is closed subgroup, we conclude that $\varphi(t)$
must actually belong to the stabilizer for all values of~$t$,  as
desired.
      \end{proof}

\section{Proof of Theorem~\ref{theorem:classification}}
\label{MainPfSect}

We assume the notation of \pref{PrelimNotation}. Recall, in
particular, that $\mu$ is carried by $\Xh$, and that the group $A
\hor$ acts on $\Xh$.

\begin{prop} \label{FibersOnHoriz}
      Almost every fiber of~$\mu$ over $\cX^d$ is supported on finitely
many orbits of $\hor \ominus X$.
      \end{prop}

\begin{proof}
     Because $\mu$ is an ergodic probability measure, it suffices to
show that almost every fiber is supported on countably many such
orbits. For $\Omega_\rho$ as in Proposition~\ref{DivStab}, we know
$\bigcup_{n=N}^\infty \Omega_{1/n}$ is conull, so it suffices to show,
for each $\rho > 0$, that each fiber of $\Omega_\rho$ is contained in
the union of countably many orbits of $\hor \ominus X$. Suppose
not. (This will lead to a contradiction.) Because any uncountable set
contains one if its accumulation points, there exist $(M',p) \in
\Omega_\rho$ and a sequence $\{p_n\}$ in~$M'$, such that
      \begin{itemize}
      \item $(M',p_n) \in \Omega_\rho$,
      \item $\pi_d(M',p_n) = \pi_d(M',p)$,
      \item $(M',p_n) \to (M',p)$,
      and
      \item $(M',p_n) \notin (\hor \ominus X) (M',p)$.
      \end{itemize}
      Because $\pi_d(M',p_n) = \pi_d(M',p)$ and $(M',p_n) \to (M',p)$,
we  may write $(M',p_n) = w_n(M',p)$ for some $w_n \in 0^d \times
(\real^2)^{k-d}$ with $w_n \to e$. By assumption, we know $w_n \notin
\hor$, so $H_1(w_n)$ is a nonzero element of $\hor \cap \bigl( 0^d
\times (\real^2)^{k-d} \bigr) = \hor \ominus X$. Because
     $$ \frac{H_1(w_n)}{\|H_1(w_n)\|} = H_{1/\|H_1(w_n)\|}(w_n) ,$$
     and  Proposition~\ref{DivStab}
implies that $H_{1/\|H_1(w_n)\|}(w_n)$ converges to an element of
$\Stab_{A\hor}(\mu)^\circ$,
     we conclude that $\Stab_{\hor \ominus X}(\mu)^\circ$ is nontrivial.
This
contradicts the definition of~$X$.
      \end{proof}

\begin{prop} \label{FinFibOverXd}
      After restricting to an appropriate conull subset~$\Omega_0$ of
$\Xh$, each fiber of~$\pi_d$ is finite.
      \end{prop}

\begin{proof}
      We know, from Proposition~\ref{FibersOnHoriz}, that almost every
fiber of~$\pi_d$ is carried by only finitely many orbits of~$\hor 
\ominus
X$. (From Theorem~\ref{KMS:UniqErg}, we may assume that each of these is
an embedded copy of $\hor \ominus X$.) Letting
   $$(\hor \ominus X)_i^+ = 0^{d+i-1} \times \real^+ \times 0^{k-d
- i} ,$$
   we may define a measurable function $\xi_i \colon \Xh \to [0,1]$
by $\xi_i(\omega) = \mu_{\pi_d(\omega)} \bigl( (\hor \ominus
X)_i^+(p) \bigr)$. This function is essentially $U$-invariant, so it
must be essentially constant. Because this is true for all~$i$, we
conclude that~$\pi_d$ is carried by a single point in each orbit of
$\hor \ominus X$. Since there are only finitely many such orbits to
consider, we conclude that a.e.\ fiber consists of a finite number of
atoms, as desired.
      \end{proof}

\begin{prop} \label{AInvt}
      We may assume $\mu$ is $A$-invariant.
      \end{prop}

\begin{proof}
   Choose $\Omega_\rho$ as in Proposition~\ref{DivStab}, with $\rho =
0.99$.
  From Corollary~\ref{UniqErgHoriz}, we know that $\mu$ projects to the
Lebesgue measure on $\cX^d$. Furthermore, by passing to a conull
subset, we may assume $\Omega_\rho$ has finite fibers over $\cX^d$
\see{FinFibOverXd}. Thus, it is easy to see that there exist
$(M,p) \in \Omega_\rho$, $\{v_n\} \subset V \smallsetminus \{e\}$,
and $\{w_n\} \subset 0^d \times (\real^2)^{k-d}$, such that $v_n
w_n(M,p) \in \Omega_\rho$, $v_n \to e$, and $w_n \to e$. Then, in the
notation of \pref{DivNotn}, with $g_n = v_n$, and choosing $s_n$
appropriately, we have $a_{n,s_n} H_{s_n}(w_n) \in A(\hor \ominus X)$
(cf.\ pf.\ of \ref{FibersOnHoriz} to see that $H_{s_n}(w_n) \in (\hor
\ominus X)$). We conclude, from Proposition~\ref{DivStab}, that the
identity component of $\Stab_{A\hor}(\mu) \cap A(\hor \ominus X)$ is
nontrivial. Because the identity component of $\Stab_{A\hor}(\mu)
\cap (\hor \ominus X)$ is trivial (by definition of~$X$), we conclude
that $\Stab_{A\hor}(\mu)$ contains a one-parameter subgroup that is
not contained in~$\hor$. Any such subgroup is conjugate to~$A$ (via an
element of~$\hor$). Thus, by replacing $\mu$ with a translate
under~$\hor$, we may assume $\mu$ is $A$-invariant.
      \end{proof}

\begin{lem}
     $\Xv$  is conull.
    \end{lem}

\begin{proof}
   By passing to a quotient, we may assume $k = 1$.
   For each nonzero vector $w \in \real^2$, let
   $$ \sing^w_g = \{\, p \in g M \mid (p + \real w) \cap g \sing \neq
\emptyset \,\} .$$
   Note that $\sing^w_{ug} = u \sing^{u^{-1} w}_g$.

   Suppose there is a subset~$E$ of positive measure in~$G$, such
that $\mu_{gM}(\sing^{(0,1)}_g) \neq 0$ for $g \in E \Gamma$. Then the
Pointwise Ergodic Theorem implies, for a.e.\ $g_0 \in G$, that we
have $u g_0 \in E \Gamma$ for all~$u$ in a non-null subset~$U_0$
of~$U$. Furthermore, because $\mu$ is $U$-invariant, we may assume
$\mu_{g_0 M} = u_* \mu_{u g_0 M}$ for all $u \in U_0$. Therefore
   $$ \mu_{g_0 M}(\sing^{u^{-1}(0,1)}_{g_0})
   = \mu_{u g_0 M}(u \sing^{u^{-1}(0,1)}_{g_0})
   = \mu_{u g_0 M}(\sing^{(0,1)}_{u g_0})
   \neq 0$$
   for all $u \in U_0$.
   This contradicts the fact that, because $\sing^{w_1}_{g_0} \cap
\sing^{w_2}_{g_0}$ is countable whenever $\real w_1 \neq \real w_2$,
we have $\mu_{gM}(\sing^w_g) = 0$ for all but countably many choices
of the line~$\real w$.
    \end{proof}

\begin{prop}[cf.\ {\cite[Cor.~8.4]{MargulisTomanov-Ratner},
\cite[Cor.~5.5.2]{Morris-Ratner}}] \label{AllInVVert}
      There is a conull subset $\Omega$ of~$\Xv$, such that
      $$(V \ver \omega) \cap \Omega = (V Y \omega) \cap \Omega ,$$
      for all $\omega \in \Omega$.
      \end{prop}

\begin{proof}
      Let $\Omega$ be a generic set for for the action of~$A$ on~$\Xv$;
thus,
$\Omega$ is conull and, for each $\omega \in \Omega$,
      $$ \mbox{$a^t \omega \in \Omega_\rho$ for most $t \in \real^+$} .$$
      Given $(M,p), (M',p') \in \Omega$, such that $(M',p') = v y (M,p)$
with
$v
\in V$ and $y \in \ver$, we wish to show $y \in Y$.

Choose a sequence $t_n \to \infty$, such that $a^{t_n}(M,p)$ and
$a^{t_n}(M',p')$ each belong to~$\Omega_\rho$. Because $t_n \to
\infty$ and  $V \ver$ is the foliation that is contracted
by~$a^{\real^+}$, we know that $a^{-t_n}(v y) a^{t_n} \to e$.
Furthermore, because $A$ acts on the Lie algebra of~$V$ with twice
the weight that it acts on the Lie algebra of~$\ver$, we see that
$\|a^{-t_n} v a^{t_n}\|/| a^{-t_n} y  a^{t_n}\| \to 0$. Thus,
letting $s$ be within a constant multiple of $1/\| a^{-t_n} y
a^{t_n}\|$, we see, in the notation of \pref{DivNotn}, with $g_n =
a^{-t_n} v a^{t_n}$ and $w_n = a^{-t_n} y a^{t_n}$, that $a_{s_n}(g_n)
\to e$, but $H_{s_n}(w_n) \not\to e$. Thus, Proposition~\ref{DivStab}
asserts that  $H_{s_n}(w_n)$ converges to a
nontrivial element of $\Stab_{\hor}(\mu)^\circ = X$. Since
$H_{s_n}(w_n) = H_{s_n}(a^{-t_n} y  a^{t_n})$ is a scalar multiple of
$H_1(y)$, we conclude that $H_1(y) \in X$. Therefore
      $(u^1 - \Id) y = H_1(y) \in X$,
      so $y \in Y$.
      \end{proof}

We require the following entropy estimate.

\begin{lem}[cf.\ {\cite[Thm.~9.7]{MargulisTomanov-Ratner},
\cite[Prop.~2.5.11]{Morris-Ratner}}] \label{EntLem}
      Suppose $W$ is a closed connected subgroup
of~$V \ver$ that is normalized by~$a \in A^+$, and let
      $$J(a^{-1}, W) = \det \bigl( (\Ad a^{-1})|_{\Lie W} \bigr)$$
      be the Jacobian of~$a^{-1}$ on~$W$.
      \begin{enumerate}

      \item \label{EntLem-Inv->}
      If $\mu$ is $W$-invariant, then $h_\mu(a) \ge \log J(a^{-1}, W)$.

      \item \label{EntLem-GminInW}
      If there is a conull, Borel subset~$\Omega$ of $\cX^k$, such
that $\Omega \cap V \ver  \omega \subset W \omega$, for
every $\omega \in \Omega$, then $h_\mu(a) \le \log J(a^{-1},
W)$.

      \item \label{EntLem-equal}
      If the hypotheses of \pref{EntLem-GminInW} are satisfied, and
equality holds in its conclusion, then $\mu$~is $W$-invariant.

      \end{enumerate}
      \end{lem}

\begin{prop}[cf.\
{\cite[Step~1 of 10.5]{MargulisTomanov-Ratner}},
{\cite[Prop.~5.6.1]{Morris-Ratner}}] \label{VYInvt}
       $\mu$ is $V Y$-invariant.
      \end{prop}

\begin{proof}
       From Lemma~\fullref{EntLem}{Inv->}, with $a^{-1}$ in the
role of~$a$, we have
      $$ \log J(a, UX) \le h_\mu(a^{-1}) .$$
      From Proposition~\ref{AllInVVert} and
Lemma~\fullref{EntLem}{GminInW}, we have
      $$ h_\mu(a) \le \log J(a^{-1}, V Y) .$$
      Combining these two inequalities with the facts that
      \begin{itemize}
      \item
$h_\mu(a) = h_\mu(a^{-1})$
      and
      \item $J(a, UX) = J(a^{-1}, VY)$,
      \end{itemize}
      we have
      $$ \log J(a, UX) \le h_\mu(a^{-1}) = h_\mu(a) \le
\log J(a^{-1}, VY) = \log J(a, UX) .$$
      Thus, we must have equality throughout, so the desired
conclusion follows from Lemma~\fullref{EntLem}{equal}.
      \end{proof}

\begin{prop}
      $\mu$ is the Lebesgue measure on a single orbit on $\cX^k_0$ of
the pseudogroup $G \ltimes \langle \Phi^k_X, \Phi^k_Y \rangle$.
       \end{prop}

\begin{proof}
      We know:
      \begin{itemize}
      \item $U$ preserves~$\mu$ (by assumption),
      \item $X$ preserves~$\mu$ (by definition),
      \item $A$ preserves~$\mu$ \see{AInvt},
      and
      \item $V Y$ preserves~$\mu$ \see{VYInvt}.
      \end{itemize}
      Therefore, $\mu$ is preserved by the pseudogroup $G \ltimes
\langle \Phi^k_X, \Phi^k_Y \rangle$ generated by these maps. Because
      \begin{itemize}
      \item this pseudogroup is transitive on the quotient $\cX^d_0$,
      and
      \item $\mu$ has finite fibers over $\cX^d_0$ \see{FinFibOverXd},
      \end{itemize}
      this implies that some orbit of the pseudogroup has positive
measure.
By
ergodicity of~$U$, then this orbit is conull.
      \end{proof}

\begin{rmk}
  To obtain the conclusions of Theorem~\ref{theorem:classification}, we
let $W = X + Y$. Then $\mu$ is supported on the $(G \ltimes
\semigrp^k_W)$-orbit of some point $(p_1,\ldots,p_k)$ in~$\cX^k$. Note
that, by choosing $\dim W$ to be minimal, we can guarantee that whenever
$p_i$ is a singular point of~$M$, the subspace $W$ projects to~$0$ in 
the
$i^\text{th}$ coordinate of~$(\real^2)^k$. Therefore, the dimension of
the orbit is equal to the dimension of the pseudosemigroup.
  \end{rmk}

\section{Countability}
\label{CountabilitySection}

For our application, we need the following analogue of
\cite[Cor.~A(2)]{RatnerEqui}.

\begin{prop}
\label{prop:countable}
The set of subspaces $W$ occurring in
Theorem~\ref{theorem:classification} is countable. For each such $W$,
the set of closed orbits of $\SL(2,\real) \ltimes \semigrp_W^k$ is
countable.
\end{prop}

\begin{lem}
\label{lemma:Ws:countable}
The set of $G$-invariant subspaces $W$ of $(\real^2)^k$ such that
there exists $p \in M^k$ with $\cO = \semigrp_W^k p$ closed is 
countable.
\end{lem}

\begin{proof} Let $2d$ be the dimension of $W$. After possibly
renumbering the factors, we may assume that
\begin{displaymath}
W \cap \left( (0,0)^d \times (\real^2)^{k-d} \right) = \emptyset.
\end{displaymath}
Then if we denote elements of $(\real^2)^k$ by $(v_1, \dots, v_k)$
where each $v_j \in \real^2$, then
$W$ is given by the following equations: for $d+1 \le j \le k$,
\begin{equation}
\label{eq:W}
v_j = \sum_{i=1}^d \alpha_{ji} v_i
\end{equation}

Recall that the linear holonomy map $\hol \colon H_1(M,\integer) \to
{\mathbb {C}} \cong
\real^2$ is given by $\hol(\gamma) = \int_\gamma \omega$, where
$\omega = dx + i dy$ is the holomorphic $1$-form that determines the
flat structure on $M$. Let $\Delta \subset \real^2$ denote the image
of $\hol$, and let $F$ denote the set of real numbers $r$ such that
there exist nonzero $v_1 \in \Delta$, $v_2 \in \Delta$ with $v_1 = r 
v_2$.
Then $F$ is clearly a countable set. We will show that each
$\alpha_{ji}$ belongs to $F \cup \{0\}$.

Let $\pi_d \colon M^k \to M^d$ denote projection onto the first $d$
factors. Note that the intersection of $\cO$ with each fiber of
$\pi_d$ is finite.

Now pick $i$, $1\le i \le d$, and $j$, $d+1 \le j \le k$. We may
assume that $\alpha_{ji} \ne 0$. Choose $p = (p_1, \dots, p_k) \in
\cO$ such that $p_i$ and $p_j$ are non-singular.  Let $\gamma$ be any
element of $H_1(M,\integer)$ with $\hol(\gamma) \ne 0$.  We represent
$\gamma$ by a piecewise linear closed curve on $M$ beginning and
ending at $p_i$ and not passing through any singularities; we will
also denote this representative by $\gamma$.  We obtain a closed curve
$\gamma_i \in M^d$ by keeping $p_m$ fixed for $1\le m \le d$, $m\ne
i$.  Because $\cO$ is a branched cover of $M^d$, $\gamma_i$ lifts to a
closed curve $\tilde{\gamma}_i$ in $\cO$. Let $\gamma'$ denote the
projection of $\tilde{\gamma}_i$ to the $j^\text{th}$ factor.

We wish to calculate $\hol(\gamma')$, so let us describe $\gamma'$
more precisely. The curve $\gamma$ is a collection of segments
connecting points $p_i = q_0, q_1, \dots, q_{n-1}, q_n = p_i$, with
$q_{m+1} = \widehat{\phi}_{w_m}(q_m)$, $w_m \in \real^2$. Then
$\gamma'$ is a collection of segments connecting the points $p_j =
q_0', q_1', \dots, q_{n-1}', q_n'$, with $q'_{m+1} =
\widehat{\phi}_{\alpha_{ji}w_m}(q_m')$. By perturbing the $w_m$, we
can make sure that $\gamma'$ is well defined and is not passing
through any singularities.

By construction, the endpoint $q_n'$ of $\gamma'$ belongs to the
finite set $\pi_d^{-1}(\pi_d(p)) \cap \cO$.  After replacing $\gamma$
by an integer multiple, we may assume that $\gamma'$ is closed. But,
in view of the explicit description of $\gamma'$, $\hol(\gamma') =
\alpha_{ji} \hol(\gamma)$, so $\alpha_{ji} \in F$.
\end{proof}

In the rest of this
section we will abuse notation by writing $p+v$ for 
$\widehat{\phi_v}(p)$.

\begin{lem}
\label{lemma:periodic:countable}
Let $M$ be a Veech surface, and let $\Gamma$ be the Veech group of
$M$. A point $p$ is called a periodic point if the $\Gamma$ orbit of
$p$ is finite. Then the set of periodic points is countable.
\end{lem}

\begin{rmk}
  When $M$ is non-arithmetic, which is the only case that we need to
discuss, it is proven in \cite{Gutkin:Hubert:Schmidt} that the number
of periodic points is countable (in fact, \emph{finite}). The
following generalization of this statement also follows from the
results of \cite{Gutkin:Hubert:Schmidt}, but we include a short
proof of as a warm up to the proof of
Proposition~\ref{prop:countable}.
   \end{rmk}

\begin{proof} It is clearly enough to show that for each $n \in
\natural$, the set $\cP_n$ of points of period $n$ is countable. To
do this it is enough to show that for each point $p \in \cP_n$,
there exists a neighborhood $U$ of $p$ that does not contain any
other points of $\cP_n$. Suppose the last statement is false. Then
there exists a sequence of points $p_j \in \cP_n$ such that $p_j
\to p$. We may assume after passing to a subsequence that the $p_j$
approach $p$ from some given direction $w$ (i.e. that $\lim
\frac{p_j-p}{\|p_j -
   p\|} = w$). Let $\Gamma'$ denote the intersection of all the index
$n$ subgroups of $\Gamma$. Then, since $\Gamma$ is finitely generated,
$\Gamma'$ is of finite index in $\Gamma$ and for each $\gamma' \in
\Gamma'$, and all $j$, $\gamma'(p_j) = p_j$. Then each element of
$\Gamma'$ must fix $w$. This contradicts the fact that $\Gamma'$,
being a finite index subgroup of $\Gamma$, is Zariski dense in $G$.
\end{proof}

\subsection*{Proof of Proposition~\ref{prop:countable}.}
It remains to prove the following assertion: Let $W \subset
(\real^2)^k$ be an $G$-invariant subspace. Then the set
$\cH$ of closed orbits of $\Gamma \ltimes \semigrp_W^k$ is countable.

We triangulate $M$, with the vertices at the singular points. This
yields a cell decomposition of $M^k$ in which the cells
$\Delta_1, \dots \Delta_m$ of maximal dimension are products of
triangles. Let $\Delta_i^0$ denote the interior of $\Delta_i$, and
let $M_0^k$ denote the union of the $\Delta_i^0$. For $p \in M^k$, let
$\delta(p)$ denote the distance between $p$ and the complement of
$M_0^k$ (i.e. the distance to the boundary of the cell containing
$p$).

Let $2d = \dim W$, and let $W^\perp$ be any $G$-invariant complement
to $W$. We may assume that $W$ is given by the equations
(\ref{eq:W}).  In view of Lemma~\ref{lemma:periodic:countable} we may
also also assume that $W$ has dense projection onto any of the
$\real^2$ factors (i.e. for a fixed $j$, not all $\alpha_{ji}$ are
$0$). Then, for any $\cO \in \cH$, $\cO \cap M_0^k$ is dense in $\cO$.

Let $n_1, \dots, n_m$ be an $m$-tuple of non-negative integers, and
let $\cH(n_1, \dots, n_m)$ denote the set of
of $\cO \in \cH$ such that  $\cO \cap \Delta_i^0$
has exactly $n_i$ connected components.

Now suppose $\cH$ is uncountable. Then there exist $n_1, \dots, n_m$
such that $\cH(n_1, \dots, n_m)$ is uncountable. Then by compactness,
there exist $\cO$ in $\cH(n_1, \dots, n_m)$ such that for every 
$\epsilon >
0$ there exists $\cO' \in \cH(n_1,\dots, n_m)$ such that the Hausdorff
distance between $\cO$ and $\cO'$ is less then $\epsilon$.
Let $\rho$ be the minimum over $i$ of the
minimal distance between connected components of $\cO \cap \Delta_i^0$.

Let $n = n_1 + \dots + n_m$, and number all the connected components of
the intersection of $\cO$ with the interiors of the cells as
$\cO_i$, $1 \le i \le n$.
Let $\gamma_1, \dots, \gamma_s$ denote the generators of $\Gamma$. We
may choose a point $p_i$ in each $\cO_i$ such that for all $j$, $1 \le
j \le s$, $\gamma_j p_i$ is in the interior of some component $\cO_l$,
where $l$ depends on $i$ and $j$.

Let $C = \max_{1 \le j \le s} \|\gamma_j\|$.
Now choose $\epsilon > 0$ so that:
\begin{itemize}
\item $ C \epsilon < \rho/3$.
\item For any $i$, $1 \le i \le m$, we have $\delta(p_i) >
   2 C \epsilon$.
\item For each $i$, $1 \le i \le m$ and each $j$, $1 \le j \le s$, we
   have $\delta(\gamma_j p_i) > 2 C \epsilon$.
\end{itemize}

Now choose $\cO' \in \cH(n_1, \dots, n_m)$ so that the Hausdorff
distance between $\cO'$ and $\cO$ is less then $\epsilon$.  Note that
if $q \in \cO_i$ with $\delta(q) \ge 2 C \epsilon$ there exists a
unique $v_i \in W^\perp$, such that $\|v_i\| \le C \epsilon$ and
$q+v_i \in \cO'$. Also $v_i$ does not depend on the choice of $q$, and
$\|v_i \| \le \epsilon$. Furthermore $v_i \ne 0$ since $\cO$ and
$\cO'$ cannot share a point in $M_0^k$.  Let $V$ denote the finite set
$\{v_1, \dots, v_n \}$.

We now claim that each generator $\gamma_j$ preserves the set
$V$. Indeed consider the points $p_i \in \cO$ and $p_i+v_i \in \cO'$.
Since both $\cO$ and $\cO'$ are $\Gamma$-invariant, we must have
$\gamma_j p_i \in \cO$ and $\gamma_j (p_i+v_i) \in \cO'$.
By construction, $\gamma_j p_i \in \cO_l$, and $\delta(\gamma_j p_i) >
2 C \epsilon$.
Recall that $v_l$ is the only vector in $W^\perp$ of norm at
most $C \epsilon$ such that $\gamma_j p_i + v_l \in \cO'$. But
$\gamma_j (p_i+v_i) = \gamma_j p_i + \gamma_j v_i \in \cO'$, and
$\|\gamma_j v_i \| \le \|\gamma_j\| \| v_i\| \le C \epsilon$. Also
$\gamma_j v_i \in W^\perp$, since $W^\perp$ is $G$-invariant.
Thus $\gamma_j v_i = v_l$.

We have proved that for each generator
$\gamma_j$, we have $\gamma_j V \subseteq V$. This immediately implies 
that
$\Gamma V = V$. Then a finite index subgroup of $\Gamma$ will fix a
single vector in $V$, which contradicts the fact that $\Gamma$ is
Zariski dense in $G$.
  \ep\medskip

\section{Averages over large circles} \label{AvgCircleSect}

Let $m_K$ denote the Haar measure on $\SO(2) \subset G$. For $x \in
\cX^k$ and
$t > 0$, let
\begin{displaymath}
\nu_t = \nu_{t,x} = a_t m_K \delta_x
\end{displaymath}
where $\delta_x$ is the atomic probability measure supported at $x$,
and $a_t = \begin{pmatrix} e^{t} & 0 \\ 0 & e^{-t} \end{pmatrix}$.
Then each $\nu_t$ is a probability measure on $\cX^k$. We can think of
$\nu_t$ as the measure supported on a circle of radius $t$ inside the
$G$-orbit through $x$. In this section we prove the following theorem:

\begin{thm}
\label{theorem:largecircles}
  Suppose $x \in \cX^k_0$.  Then there exists a $G$-invariant
subspace~$W$ of $(\real^2)^k$ such that
  \begin{enumerate}
  \item the $G \ltimes
\semigrp^k_W$ orbit through $x$ is closed,
  and
  \item $\lim_{t \to \infty} \nu_t = \mu$, where $\mu$ is
  Lebesgue measure on this orbit.
  \end{enumerate}
\end{thm}

\begin{rmk}
\label{remark:largecircles}
If $W = (\real^2)^k$, then $\mu$ is the Lebesgue measure on~$\cX^k$.
  \end{rmk}

\begin{lem}[{Invariance under a unipotent}]
\label{lemma:circles:unipotent}
  Suppose $t_i \to \infty$. Then there is a subsequence $t_{i_j}$ such
that the measures $\nu_{t_{i_j}}$ converge to a probability measure
$\nu_\infty$ that is invariant under the unipotent element $u =
\begin{pmatrix} 1 & 1 \\ 0 & 1 \end{pmatrix}$ of $G$.
  \end{lem}

\begin{proof}
  It follows from \cite[Corollary~5.3]{EM} that there is a
subsequence $t_{i_j}$ such that the measures $\nu_{t_{i_j}}$
converge to a probability measure $\nu_\infty$. We can find
$\theta_j \to 0$ such that $a_{t_{i_j}} r_{\theta_j}
a_{t_{i_j}}^{-1}$ converges to $u$. (Recall that $r_\theta$ is the
$2 \times 2$ matrix representing rotation by $\theta$). Now the
measures $\nu_{t_{i_j}} = (a_{t_{i_j}} \nu_K a_{t_{i_j}}^{-1} )
a_{t_{i_j}} \delta_x$ are $a_{t_{i_j}} r_{\theta_j}
a_{t_{i_j}}^{-1}$ invariant, hence $\nu_\infty$ is $u$ invariant.
  \end{proof}

\begin{assump}
   Assume $\nu_\infty$ is not the Lebesgue measure on~$\cX^k$.
  \end{assump}

\subsection*{Application of the measure classification theorem.}
  Note that we do not know at this point whether $\nu_\infty$ is ergodic.
However, standard results (using $u$-invariance) imply that
$\nu_\infty$ projects to Lebesgue measure in $G/\Gamma$.

\begin{notation}
  For convenience, if $B \subset (\real^2)^k$ and $X \subset \cX^k$, let
  $$ B^{-1} X = \bigcup_{v \in B} (\withsing{\phi}^k_v)^{-1}(X) .$$
  \end{notation}

By Theorem~\ref{theorem:classification}, and by
Proposition~\ref{prop:countable}, there exists a $G$-invariant proper
subspace~$W$ of~$(\real^2)^k$ and an orbit $\cO$ of $G \ltimes
\semigrp^k_W$, such that
  $$\nu_\infty( \hor^{-1} \cO) > 0 .$$
  We will show that this implies that $x \in \cO$. In that case, the
entire $G$-orbit of $x$ lies in $\cO$, so $\nu_\infty(\cO) = 1$.
Furthermore, we show that as long as $W$ was chosen as small as
possible, $\nu_\infty$ must be Lebesgue measure on $\cO$.

\subsection*{Projection and fiber measures.}
  We choose $W$ to be of minimal dimension.
  From the structure of the $G$-invariant subspaces on $(\real^2)^k$, we
see that $\dim W = 2 d$, $0\le d < k$, and after renumbering the
factors, we can make sure that
  $$ \mbox{$W$ projects surjectively to $(\real^2)^d \times 0^{k-d}$.}
$$
  Thus, $(0,0)^d \times (\real^2)^{k-d}$ is complementary to $W$, and
$0^d \times \real^{k-d}$ is complementary to $W \cap \hor$ in~$\hor$.

\begin{lem}
\label{lemma:support:in:interval}
  There exists
$\epsilon > 0$ and a box
$$B = \{0\}^d \times [\alpha_{d+1}, \beta_{d+1}] \times \dots \times
[\alpha_k, \beta_k] \subset \hor $$
such that $\nu_\infty( B^{-1} \cO ) > 2 \epsilon$.
  \end{lem}

\begin{proof}
  Let $\cO_0$ be the (unique) orbit of $G \ltimes \psgrp^k_W$ that is
open and dense in~$\cO$. (In other words, $\cO_0$ consists of the
elements in~$\cO$ of which as few coordinates as possible are singular
points.) Note that
  \begin{equation} \label{lemma:support:in:intervalPf-Winvt}
  W^{-1} \cO_0 = \cO_0 .
  \end{equation}
  By the minimality of $\dim W$, we see that $\nu_\infty \bigl(
\hor^{-1} (\cO \smallsetminus \cO_0) \bigr) = 0$. Hence $\nu_\infty (
\hor^{-1} \cO_0) > 0$. By combining this with
\pref{lemma:support:in:intervalPf-Winvt} and the fact that $\hor = W +
(0^d \times \real^{k-d})$, we conclude that there is a box $B \subset
0^d \times \real^{k-d}$, such that $\nu_\infty( B^{-1} \cO_0 ) > 0$.
Since $\cO_0 \subset \cO$, then $\nu_\infty( B^{-1} \cO ) > 0$, as
desired.
  \end{proof}

As in the previous sections, let $\pi_d \colon \cX^k \to \cX^d$ be the
natural projection onto the first $d$ coordinates. For $z \in \cX^d$, we
let $F_z = \pi_d^{-1}(z) \cap \cO$. Note that $F_z$ is a finite set.

We claim that
  \begin{equation} \label{ProjToLebOnXd}
  \mbox{$\nu_\infty$ projects to the Lebesgue measure on~$\cX^d$.}
  \end{equation}
  To see this, note that, because $W$ is a proper subspace of
$(\real^2)^k$, we have $d < k$. Hence, by induction on~$k$, we may
assume there is a  $G$-invariant subspace~$W_d$ of~$(\real^2)^d$, such
that the projection of~$\nu_\infty$ to~$\cX^d$ is the Lebesgue measure 
on
the $G \ltimes \semigrp^d_{W_d}$ orbit~$\cO_d$ through $\pi_d(x)$. Then
$\pi_d^{-1} ( \hor^{-1} \cO_d)$ is conull for~$\nu_\infty$, so
  $$ \nu_\infty \bigl( \pi_d^{-1} ( \hor^{-1} \cO_d) \cap \cO \bigr)
  = \nu_\infty (\cO)
  \neq 0 .$$
  From the minimality of $\dim W$, we conclude that $W_d = (\real^2)^d$.
Therefore $\cO_d = \cX^d$, which establishes the claim.

\begin{assump} \label{AssumeXNotinO}
  We may assume $x \notin \cO$. (Otherwise, from the fact that $\cO$ is a
branched cover of~$\cX^d$ (and Lemma~\ref{FinCovClassify} below), we
would immediately conclude that $\nu_\infty$ is the Lebesgue measure
on~$\cO$, as desired.) This will lead to a contradiction.
  \end{assump}

\subsection*{The key estimate.}
For $L_1 > 0$, $\delta_1 > 0$ let
\begin{displaymath}
B(\delta_1, L_1) = \left\{
       (x,y) \in (\real^2)^k  \left| \
       \begin{aligned}
       x_i = y_i = 0, \qquad & \text{ for $1 \le i \le d$}, \\
       |x_i| \le L_1 \text{ and }
       |y_i| \le \delta_1, & \text{ for $d+1 \le i \le k$}.
      \end{aligned}
\right.\right\}
\end{displaymath}

\begin{lem}[{The key estimate}]
\label{lemma:trig:polynomial}
Suppose $B \subset B(\delta_1, L_1) \subset B(\delta, L)$,
where $B$ is as defined in Lemma~\ref{lemma:support:in:interval}.
Suppose also that $\rho > 0$, $\epsilon < 1$, $\delta_1 < \epsilon
\delta/5$, and
$L_1 < \epsilon L /5$. Then there exists $t_0$ depending only on
$\rho, \delta, L$ such that for any $t > t_0$ and any $v
\in B(\delta_1,L_1)$ with
\begin{equation}
\label{eq:far:from:hor}
d(v,\hor) > \frac{k}{5}e^{-t} \rho,
\end{equation}
we have
\begin{displaymath}
\bigl| \{\,\theta \mid a_t r_\theta  a_t^{-1} v \in B(\delta_1, L_1)\,\}
\bigr| \le \frac{\epsilon}{2} \bigl| \{\, \theta \mid  a_t r_\theta
a_t^{-1} v \in B(\delta, L) \,\}\bigr|,
\end{displaymath}
where $r_\theta = \begin{pmatrix} \cos \theta & \sin \theta \\ -\sin
     \theta & \cos \theta \end{pmatrix}$.
\end{lem}

\begin{figure}
  \begin{center}
  \input{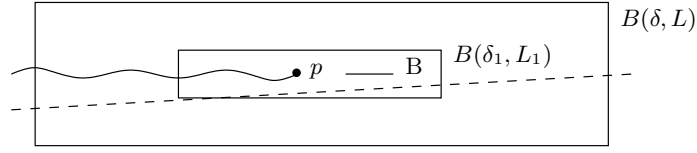}
  \caption{The time the ellipse (drawn here as a dotted line) spends
inside the small box $B(\delta_1, L_1)$ is at most $\epsilon$ times the
time the ellipse spends in the larger box $B(\delta,L)$. In
Lemma~\ref{lemma:trig:polynomial}, this is proved as a result
in~$(\real^2)^k$. Because of Lemma~\ref{lemma:neighborhoods}, it can be
transferred to~$\cX^k$, even if the ellipse crosses the branch cut
starting at the possibly singular point~$p$.} \label{fig:boxes}
  \end{center}
  \end{figure}

\begin{proof}
If we write $v = (v_1, \dots, v_k)$, with $v_j \in \real^2$,
and also write $v_j = \begin{pmatrix} x_j \\ y_j \end{pmatrix}$ then the
condition (\ref{eq:far:from:hor}) implies that there exists at least one
$j$, $m+1 \le j \le k$ with $|y_j| > \frac{1}{5} e^{-t} \rho$. The rest
of the argument will take place in the $j^\text{th}$ factor (See
Figure~\ref{fig:boxes}).

We note that the components of the map $\theta \to a_t r_\theta
a_t^{-1} v$
are trigonometric polynomials of degree $1$. In other words, the path
$\theta \to a_t r_\theta a_t^{-1} v_j$ parametrizes an ellipse.
Let $t_0 = \max(\log \frac{5L}{\rho},0)$. Then if $t > t_0$ and $\theta
=
\pi/2$ then $a_t r_\theta a_t^{-1} v_j = \begin{pmatrix} e^{2t} y_j \\
     -e^{-2t} x_j \end{pmatrix} \not\in B(\delta, L)$. Thus, the ellipse
$\theta \to a_t r_\theta a_t^{-1} v_j$ leaves $B(\delta, L)$. Then in
view of the dimensions of the boxes, the portion of the ellipse in
$B(\delta_1, L_1)$ is at most $\epsilon/2$ times the portion of the
ellipse in $B(\delta, L)$.
\end{proof}

\begin{lem} \label{lemma:neighborhoods}
  For any $L > 0$, there exists $\delta > 0$ and a compact subset $E$ of
$\cX^d$ with $\nu_\infty \bigl( \pi_d^{-1}(E) \bigr) > 1- \epsilon/4$,
such that
  \begin{enumerate}
  \item $B(\delta,L)^{-1} \bigl( \cO  \cap \pi_d^{-1}(E) \bigr)$  does
not contain any singular points, other than perhaps points in~$\cO$,
  and
  \item \label{lemma:neighborhoods-unique}
  for each $p \in B(\delta,L)^{-1}
\bigl( \cO  \cap \pi_d^{-1}(E) \bigr)$, there is a unique $b \in
B(\delta, L)$, such that $\withsing{\phi}^k_v(p) \in \cO$.
  \end{enumerate}
  \end{lem}

\begin{proof}
  For $z \in \cX^d$, let $S_z$ be the surface corresponding to~$z$ (so we
may write $\pi_d^{-1}(z)$ as $S_z^{k-d}$). Let $\sing_z$ denote the
singular set of $S_z$. For any $L > 0$ on the fixed surface $S_z$, there
exist only finitely many horizontal trajectories of length at most $L$
connecting points of $F_z \cup \sing_z$ to points of $F_z \cup \sing_z$.
Therefore we can find a large compact subset $E$ of $\cX^d$ such that
for any $z \in E$, $S_z$ has no horizontal trajectories of length at
most $2L$ connecting points of $F_z \cup \sing_z$ to other points of
$F_z \cup \sing_z$. Since $\nu_\infty$ projects to Lebesgue measure on
$\cX^d$ \see{ProjToLebOnXd}, we may choose $E$ so that $\nu_\infty 
\bigl(
\pi_d^{-1}(E) \bigr) > 1 - \epsilon/4$. Now we can choose $\delta > 0$
by compactness.

  Note that, because $B(\delta,L) \subset (0,0)^d \times 
(\real^2)^{k-d}$,
we have $B(\delta,L)^{-1} F_z \subset \pi^{-1}(z)$. Therefore
$B(\delta,L)^{-1} F_z \cap \cO = F_z$.
  \end{proof}

\subsection*{Completion of the proof of
Theorem~\ref{theorem:largecircles}.}
  Because $x \not\in \cO$ \see{AssumeXNotinO}, we may choose $\rho > 0$
so that $d(x, \cO) > k \rho$. We may also assume that on the surface
corresponding to~$x$, the distance between any two singular points is at
least $k \rho$. Let $B$, $[\alpha_i, \beta_i]$ and $\epsilon$ be as in
Lemma~\ref{lemma:support:in:interval}.  Choose $L_1$ so that for all
$d+1 \le i \le k$, we have $[\alpha_i, \beta_i] \subset [-L_1, L_1]$. 
Let
$L = 10 L_1/\epsilon$. Now choose $E \subset \cX^d$ and $\delta > 0$ so
that Lemma~\ref{lemma:neighborhoods} holds. Finally, choose $\delta_1  =
\epsilon \delta/10$. Assume $t > \log (5L/\rho)$. We will abuse notation
by writing $p+v$ for $\withsing{\phi}^k_v(p)$.

We claim that if $a_t r_{\theta_0} x + v \in \cO$, with $v \in
B(\delta_1, L_1)$, then (\ref{eq:far:from:hor}) holds. Indeed, we have
then
  $$ r_{\theta_0} x + a_t^{-1} v = a_t^{-1} (a_t r_{\theta_0} x + v) \in
a_t^{-1} \cO = \cO,$$
  so
  $$ |a_t^{-1} v | \ge d(r_{\theta_0}x, \cO) = d(x, \cO) > k \rho .$$
  Also,
  $$ |a_t^{-1} v | \le e^t \cdot d(v,\hor) + \frac{L_1}{e^t} .$$
  Therefore
  $$ d(v,\hor) \ge e^{-t} \left( |a_t^{-1} v | - \frac{L_1}{e^t} \right)
  > e^{-t} \left( k \rho - \frac{L_1}{5 L/\rho} \right)
  > e^{-t} \frac{k}{5} \rho .$$

Now let
\begin{displaymath}
R = \bigset{\theta }{a_t r_\theta x \in B(\delta_1, L_1)^{-1}
  \bigl( \cO  \cap \pi_d^{-1}(E) \bigr) } .
\end{displaymath}
Suppose $\theta \in R$. Let $v$ be the unique element of $B(\delta_1,
L_1)$ with $a_t r_\theta x + v \in \cO$, and let
  \begin{displaymath}
  I_\theta' = \{\, \theta' \mid a_t r_{\theta'}
r_{\theta}^{-1} a_t^{-1} v \in B(\delta, L) \,\} .
  \end{displaymath} Note that $\theta \in I_\theta'$, so we may let
$I_\theta$ be the component of~$I_\theta'$ that contains~$\theta$. By
(the proof of) Lemma~\ref{lemma:trig:polynomial}, $|I_\theta \cap R| \le
(\epsilon/2) |I_\theta|$.

We claim that if $I_{\theta_1} \neq
I_{\theta_2}$, then $I_{\theta_1} \cap I_{\theta_2}$ is disjoint
from~$R$.
  To see this, note that if $\theta' \in I_{\theta_1} \cap I_{\theta_2}$,
then there exist $v_1,v_2 \in B(\delta_1,
L_1)$, such that, letting
  $$ v_i' = a_t r_{\theta'} r_{\theta_i}^{-1} a_t^{-1} v_i ,$$
  we have $v_i' \in B(\delta,L) $ and
  $$ a_t r_{\theta'} x + v_i'
  = a_t r_{\theta'} r_{\theta_i}^{-1} a_t^{-1} (a_t r_{\theta_i} x + v_i)
  \in a_t r_{\theta'} r_{\theta_i}^{-1} a_t^{-1} \cO
  = \cO .$$
  Now if $I_{\theta_1} \neq
I_{\theta_2}$, then
  $r_{\theta}^{-1} a_t^{-1} v_1 \neq r_{\theta}^{-1} a_t^{-1} v_2$,
  so $v_1' \neq v_2'$. Lemma~\fullref{lemma:neighborhoods}{unique}
therefore implies that $a_t r_{\theta'} x \notin B(\delta,L)^{-1}
\bigl( \cO  \cap \pi_d^{-1}(E) \bigr)$, so $\theta' \notin R$.

  Since each point of $R \cap I_\theta$ is contained in a unique 
interval,
the circle is covered at most twice by the intervals $I_\theta$. It
follows that $|R| < \epsilon$. Equivalently, this means that
  $$\nu_t \Bigl( \bigl( B(\delta_1, L_1)^{-1} \cO \bigr) \cap
\pi_d^{-1}(E) \Bigr) < \epsilon.$$
  Since this holds for all sufficiently large $t$, we get
  $$\nu_\infty\Bigl( \bigl( B(\delta_1, L_1)^{-1} \cO \bigr) \cap
\pi_d^{-1}(E) \Bigr) \le  \epsilon.$$
  Since $\nu_\infty$ projects to Lebesgue measure, we know that
$\nu_\infty \bigl( \pi_d^{-1}(E) \bigr) > 1 - \epsilon/4$. Hence
$\nu_\infty \bigl( B(\delta_1, L_1)^{-1} \cO \bigr) < 5 \epsilon/4$.
This contradicts Lemma~\ref{lemma:support:in:interval}.
  \ep

  \begin{cor} \label{corollary:GOrbit}
  Suppose $x \in \cX^k_0$.  Then there exists a $G$-invariant
subspace~$W$ of $(\real^2)^k$, such that the closure of $G x$ is $(G
\ltimes \semigrp^k_W)(x)$.
      \end{cor}

\begin{proof}
  Let $W$ be as in the conclusion of Theorem~\ref{theorem:largecircles}.
Because $(G \ltimes \semigrp^k_W)(x)$ is closed and $G$-invariant, it
contains the closure of $G x$. On the other hand, the support of $\nu_t$
is a subset of $G x$, so $G x$ is dense in the support of $\lim_{t
\to \infty} \nu_t$; that is, $G x$ is dense in $(G \ltimes
\semigrp^k_W)(x)$.
  \end{proof}

Corollary~\ref{corollary:GammaOrbit} (stated at the end
of~\S\ref{MeasClassSect}) follows from \pref{corollary:GOrbit} by a
standard argument (inducing the action of $\Gamma$ to an action of~$G$).

\section{Application to counting}
\label{Application}

We now give the general setup for the counting problems we are
considering.
For additional background and more detailed definitions, see the
introduction to \cite{EMZ}.

\begin{notation}
\label{notation:counting}
\begin{itemize}

\item Let $S$ be a translation surface. A {\em saddle connection} on 
$S$ is a straight
line segment connecting two singularities. Since a saddle connection
has a well defined length and direction, each saddle connection is
associated with a non-zero vector in $\real^2$. Let $V_{sc}(S)
\subset \real^2$ denote the set of vectors in $\real^2$ that are
associated to saddle connections in $S$.

\item By a {\em regular} closed geodesic on $S$, we mean a closed 
geodesic
that does not pass through singularities.

\item
As mentioned in the introduction, any regular closed geodesic is part
of a family of freely homotopic parallel closed geodesics of the same
length. Such a family is called a {\em
   cylinder}. All the geodesics comprising a cylinder have the same
length and direction; thus we can associate to a cylinder a non-zero
vector in $\real^2$.  Note that each boundary component of a cylinder
is a union of  saddle connections. Let $V_{cyl}(S) \subset \real^2$
denote the set (with multiplicity) of vectors in $\real^2$ that are
associated to cylinders in $S$. In particular, if $S$ is a standard
torus, then $V_{cyl}(S)$ is the set of primitive vectors in
$\integer^2$.

\item For any $T >0$, let $B(T)$ denote the ball in~$\real^2$ of
radius~$T$ centered at~$0$.

\item Let $V(S)$ be a subset of $\real^2 - (0,0)$ with multiplicity; 
i.e. a
set of vectors with positive weights. The weights are usually positive
integers (e.g., we may consider saddle connections with multiplicity),
but need not be (e.g., we may weight each cylinder by the reciprocal
of its area).

\item  Let $N_V(S,T)$ denote the cardinality (with weights)
of $V(S) \cap B(T)$. We are interested in the asymptotics of
$N_V(S,T)$ as $T \to \infty$. If $V(S) = V_{sc}(S)$, we will denote
$N_V(S,T)$ by $N_{sc}(S,T)$, and if $V(S) = V_{cyl}(S)$ then, as in
the introduction, we will denote $N_V(S,T)$ simply by $N(S,T)$.

\item Recall from the introduction
that $\cH(\beta)$ denotes a {\em stratum} of translation surfaces.

\item  Let
$\cH_1(\beta)$ denote the subset of $\cH(\beta)$ consisting of the
surfaces of area $1$ (where area is taken using the associated 
translation
metric).

\item As in \S\ref{AvgCircleSect}, let
$m_K$ denote the Haar measure on $\SO(2) \subset \SL(2,\real)$.

\item
For $S \in
\cH_1(\beta)$ and
$t > 0$, let
\begin{displaymath}
\nu_{t,S} = a_t m_K \delta_S
\end{displaymath}
where $\delta_S$ is the atomic probability measure supported at $S$,
and $a_t = \begin{pmatrix} e^{t} & 0 \\ 0 & e^{-t} \end{pmatrix}$.
Then $\nu_{t,S}$ is a probability measure on $\cH_1(\beta)$.

\item Finally, for a bounded compactly supported function $f \colon
\real^2 \to \real$, let
\begin{displaymath}
\hat{f}_V(S) = \sum_{v \in V(S)} f(v) .
\end{displaymath}
The function $\hat{f}_V$ is called the {\em Siegel-Veech} transform of
$f$.

\end{itemize}
\end{notation}

\subsection*{The general counting problem.}

We now summarize the relevant results from \cite{Veech:Siegel},
\cite{EM} and \cite{EMS} that will be used in \S\ref{tribillSect}.

\begin{thm}
\label{theorem:counting:EMS:EM}
Let $S \in \cH_1(\beta)$ be a translation surface, and
suppose the following hold {\rm(}using
Notation~\ref{notation:counting}{\rm)}:
  \begin{enumerate}
  \renewcommand{\theenumi}{\Alph{enumi}}
  \renewcommand{\labelenumi}{\rm(\theenumi)}
\item \label{(A)} $V(\cdot)$ varies linearly under the $\SL(2,\real)$
   action; i.e., for all $g \in \SL(2,\real)$ and all $S \in
   \cH_1(\beta)$, we have $V(g S) = g V(S)$.

\item \label{(B)} There exists a constant $C$, such that for all $S
\in
   \cH_1(\beta)$, we have $N_V(S,2) \le C N_{sc}(S,2)$.

\item \label{(C)} As $t \to \infty$, the measures $\nu_{t,S}$ converge
to
   an $\SL(2,\real)$-invariant {\rm(}probability{\rm)} measure $\mu$.

\item \label{(D)} Let $h \colon \real^2 \to \real$ denote the
characteristic
   function of the trapezoid whose vertices are at $(1,1)$, $(0,1)$,
   $(0,1/2)$ and $(1/2,1/2)$. Let $\widetilde{\cO}$ denote the closure 
of the
   $\SL(2,\real)$ orbit of $S$. Then for any $\epsilon > 0$ and any
   compact subset $K$ of $\cH_1(\beta)$, there
   exist continuous functions $\phi_{+} \colon \widetilde{\cO} \to
\real$ and $\phi_{-} \colon \widetilde{\cO}
   \to \real$ such that for all $S \in \widetilde{\cO} \cap K$, we have
$$ \text{$\phi_-(S) \le \hat{h}_V(S)
   \le \phi_{+}(S)$ and $\int_{\cH_1(\beta)} (\phi_+ - \phi_-) \, d\mu
  < \epsilon$. }$$
\end{enumerate}
Then, the following hold:
\begin{enumerate}
  \renewcommand{\theenumi}{\roman{enumi}}
  \renewcommand{\labelenumi}{\rm(\theenumi)}
\item \label{(i)} There exists a constant $c = c(S,V)$, such that as
$T
   \to \infty$,
$$N_V(S, T) \sim \pi c  T^2.$$

\item \label{(ii)}
  We have the Siegel-Veech formula: there exists a
   constant $c$ such that for any continuous
   compactly supported $f \colon \real^2 \to \real$,
\begin{equation}
\label{eq:Siegel:Veech:general}
\int_{\cH_1(\beta)} \hat{f}_V \, d\mu = c \int_{\real^2} f .
\end{equation}
\item \label{(iii)} The constant $c$ in \pref{(i)} is the same as
the constant
   $c$ in \pref{(ii)}.
\end{enumerate}
\end{thm}

\begin{rmk}
Conclusion~\pref{(ii)} depends only on assumptions~\pref{(A)} and (some
version of)~\pref{(B)}. It was proved by W. Veech in
\cite{Veech:Siegel}, where this approach to counting on translation
surfaces was originated. The proof is reproduced in \cite[Theorem
2.2]{EM}. \end{rmk}

\begin{rmk} Assumption \pref{(B)} may be replaced by:
\begin{enumerate}
  \renewcommand{\theenumi}{\ref{(B)}$'$}
  \renewcommand{\labelenumi}{\rm(\theenumi)}
\item \label{(B')}
  There exist  constants $C > 0$ and $0 < s < 2$
   such that for all $S \in \cH_1(\beta)$, $N_V(S,2) \le
   {C}/{\ell(S)^s}$, where $\ell(S)$ is the length of the shortest
   saddle connection on $S$.
\end{enumerate}
In fact, \pref{(B')} is used in the proof of
Theorem~\ref{theorem:counting:EMS:EM} instead of \pref{(B)}. The 
assertion
that \pref{(B)} implies \pref{(B')} follows from \cite[Theorem 5.1]{EM}.
\end{rmk}

\begin{rmk}
It follows from \pref{(B')} and \cite[Theorem 5.2]{EM} that any limit 
measure of the
probability measures $\nu_{t,S}$ must be a probability measure
(see \cite[Corollary~5.3]{EM}).  Thus,
the measure $\mu$ of \pref{(C)} is automatically a probability measure.
\end{rmk}

\begin{rmk}
The assertion \pref{(D)} is a technical assumption needed since the
Siegel-Veech transform $\hat{f}$ may not be continuous even if $f$
is.
\end{rmk}

\subsection*{Outline of proof of
   Theorem~\ref{theorem:counting:EMS:EM}.}
Let $h$ be the characteristic function of the trapezoid as in 
\pref{(D)}.
We have the following lemma from calculus
(cf. \cite[Lemma 3.4]{EM}): for any $v \in \real^2$,
\begin{equation}
\label{eq:calculus}
\int_0^{2 \pi} h( a_t r_\theta v ) \, d\theta \approx \begin{cases} 
e^{-2 t} & \text{
     if $e^{t}/2 \le \|v \| \le e^t$ } \\
    0  & \text{ otherwise} .
\end{cases}
\end{equation}
If we multiply both sides of (\ref{eq:calculus}) by $T^2 = e^{2 t}$
and sum over all $v \in V(S)$, we get, under assumption \pref{(A)},
\begin{displaymath}
T^2 \int_0^{2\pi} \hat{h}_V(a_t r_\theta S) \, d\theta \approx N_V(S,T) 
- N_V(S,T/2)
  \end{displaymath}
or equivalently,
\begin{equation}
\label{eq:count:integral}
N_V(S,T) - N_V(S,T/2) \approx 2 \pi T^2 \int_{\cH_1(\beta)} \hat{h}_V 
\, d\nu_{t,S}
  . \end{equation}
(The fact that we only have approximate equality and not equality in
(\ref{eq:count:integral}) does not affect the asymptotics. See
\cite[\S{3}]{EM} for the details).

The assumption \pref{(C)} means that
for any bounded continuous function $\phi$ on $\cH_1(\beta)$,
\begin{equation}
\label{eq:circles:bounded:func}
\lim_{t\to \infty}
\int_{\cH_1(\beta)} \phi \, d \nu_{t,S} = \int_{\cH_1(\beta)} \phi \, 
d\mu
  . \end{equation}
We would like to apply (\ref{eq:circles:bounded:func}) to $\hat{h}_V$,
which is neither bounded nor continuous. The fact that $\hat{h}_V$ is
not continuous is handled by assumption \pref{(D)}. To handle the fact 
that
$\hat{h}_V$ is not bounded, we decompose $\hat{h}_V = h_1 + h_2$, where
$h_1$ is bounded and $h_2$ is supported outside of a large compact set. 
Then
the contribution of $h_2$  can be shown to be negligible using 
\cite[Theorem
5.2]{EM}, in view of assumption \pref{(B')}.  The details of this 
argument
are given in \cite[\S 2]{EMS}.

Now applying (\ref{eq:circles:bounded:func}) with $\phi = \hat{h}_V$
and substituting into (\ref{eq:count:integral}), we get
\begin{displaymath}
\lim_{T \to \infty} \frac{N(S,T) - N(S,T/2)}{T^2} = 2
\pi  \int_{\cH_1(\beta)} \hat{h}_V \, d \mu .
\end{displaymath}
By iterating (replacing $T$ with $T/2, T/4, T/8, \ldots$),
and summing the resulting geometric series, we get
\begin{equation}
\label{eq:temp:asymp}
  \lim_{T \to \infty} \frac{N(S,T)}{T^2} = \frac{8 \pi}{3}
\int_{\cH_1(\beta)} \hat{h}_V \, d \mu .
\end{equation}
This implies \pref{(i)}. Now by (\ref{eq:Siegel:Veech:general}),
\begin{displaymath}
\int_{\cH_1(\beta)} \hat{h}_V \, d \mu = c \int_{\real^2} h = 
\frac{3c}{8}.
\end{displaymath}
This, together with (\ref{eq:temp:asymp}), implies \pref{(iii)}.
{\ep\medbreak}

As a corollary of Theorem~\ref{theorem:counting:EMS:EM} and
Theorem~\ref{theorem:largecircles} we have the following:

\begin{thm}
\label{theorem:branched:cover:count}
Suppose $S$ is a branched cover of a Veech surface $M$. Let $N(S,T)$
denote the number of cylinders of periodic trajectories in $S$ of
length at most $T$. Then there exists a constant $c = c(S)$ such that
as $T \to \infty$,
\begin{equation}
\label{eq:branch:asymp}
N(S,T) \sim c T^2.
\end{equation}
\end{thm}

\begin{proof}
We use Theorem~\ref{theorem:counting:EMS:EM}, with $V(\cdot) =
V_{cyl}(\cdot)$. Assumption \pref{(A)} clearly holds, and
\pref{(B)} also holds since the boundary of every cylinder contains a 
saddle
connection.

Now let $\cM$ be the
connected component of $\cM_D(\beta)$ that contains~$S$ (where
$\cM_D(\beta)$ is as in the introduction). Since $S \in \cM$ and $\cM$ 
is
closed and $\SL(2,\real)$-invariant, the support of any of the measures
$\nu_{t,S}$ is contained in $\cM$.
Also, since $\cM$ is a branched cover of the
space $\cX^k$, a measure classification theorem on $\cX^k$ automatically
yields a measure classification theorem on $\cM$ (see
Lemma~\ref{FinCovClassify} below). Thus Assumption \pref{(C)} of
Theorem~\ref{theorem:counting:EMS:EM} follows from
Theorem~\ref{theorem:largecircles}.

Finally, in our setting \pref{(D)} is automatically satisfied, since
the orbit closure $\widetilde{\cO}$ is a proper submanifold of 
$\cH_1(\beta)$,
the measure $\mu$ is Lebesgue measure on $\widetilde{\cO}$, and (after
intersecting with any compact set) the set of
discontinuities of $\hat{h}_V$ is contained in a finite union of 
submanifolds
of positive codimension in $\widetilde{\cO}$.
Thus Theorem~\ref{theorem:branched:cover:count} follows from
\pref{(i)} of Theorem~\ref{theorem:counting:EMS:EM}.
\end{proof}

\begin{lem} \label{FinCovClassify}
  Suppose $W$ is a $G$-invariant subspace of~$(\real^2)^k$, $\cO$~is
a closed orbit of $G \ltimes \semigrp^k_W$ in~$\cX^k$, and
$\widetilde{\cO}$ is a {\rm(}connected{\rm)} branched cover of~$\cO$,
such that the action of~$G$ on~$\cO$ lifts to~$\widetilde{\cO}$.

If $\nu$ is any $u$-invariant probability measure on~$\widetilde{\cO}$
that projects to the Lebesgue measure on~$\cO$, then $\nu$ is the
Lebesgue measure on~$\widetilde{\cO}$.
  \end{lem}

\begin{proof}
  Let $\mu$ and $\widetilde{\mu}$ be the Lebesgue measures on
$\cO$ and $\widetilde{\cO}$, respectively. Then, because it projects
to~$\mu$, the measure $\nu$ must be absolutely continuous with respect
to~$\widetilde{\mu}$; thus, we may write $\nu = f \widetilde{\mu}$, for
some Borel function~$f$ on~$\widetilde{\cO}$.

It is not difficult to see that $\widetilde{\mu}$ is ergodic for~$G$, so
(by decay of matrix coefficients \cite[Thm.~2.4.2, p.~29]{ZimmerBook},
or by the the Mautner phenomenon \cite[Thm.~2.2.15, p.~21]{ZimmerBook})
it is ergodic for~$u$. This implies that $f$~is constant. So $\nu =
\widetilde{\mu}$ (up to a normalizing scalar multiple).
  \end{proof}

\section{Triangular Billiards.} \label{tribillSect}
Let $n \ge 5$ be an odd integer.  As in the introduction, let
  $$ \text{$P_n$
denote the triangle with angles $\displaystyle \frac{(n-2)\pi}{2n},
\frac{(n-2)\pi}{2n}, \frac{2\pi}{n}$} $$
  and let $S_n$ denote the corresponding translation surface.  In the
rest of this section, we complete the proof of
Theorem~\ref{theorem:triangle:count} by computing the constant $c$ in
Theorem~\ref{theorem:branched:cover:count} for the case of the surface
$S_n$. Our general strategy is to use \pref{(ii)} and \pref{(iii)} of
Theorem~\ref{theorem:counting:EMS:EM}.
To pass from~$S_n$ to~$P_n$, note that $N(P_n,T) = N(S_n,T)$,
and since $S_n$ consists of $4 n$ triangles, $\area(S_n) = 4n
\area(P_n)$.

\begin{figure}
  \begin{center}
  \input{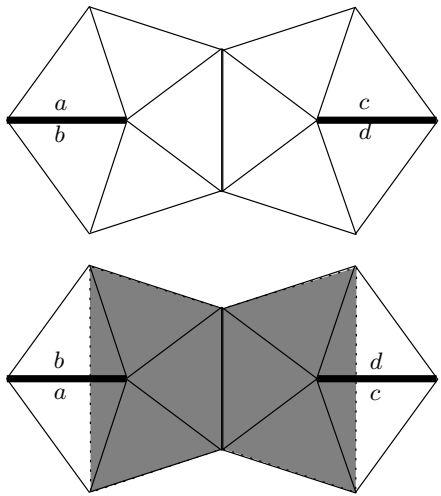}
  \caption{We draw the surface $S_n$ (for $n = 5$), tessellated by
(reflections of) the triangle $P_n$. In each of the double $n$-gon
shapes, the opposite parallel sides are identified. The bottom double
$n$-gon can be identified with the surface $X_n$. The covering map
from $S_n$ to $X_n$ is specified by the two slits (drawn as thick
lines), with identifications as shown.
  For $n = 5$, the shaded region in the bottom double pentagon is one
of the cylinders in the vertical cylinder decomposition for $X_n$; the
unshaded region in the bottom double pentagon is the other cylinder.}
  \label{fig:pentagons}
  \end{center}
  \end{figure}

The surface $S_n$ can be drawn as in Figure~\ref{fig:pentagons}. As
shown in \cite{HS1} and as one can see from the figure, $S_n$ is a
double cover of a surface
$X_n$ consisting of a double $n$-gon with opposite sides
identified. The surface $X_n$ is a Veech surface (see
\cite{Veech:eisenstein}), but $S_n$ is not (see \cite{HS1}).

\section*{The Veech surface.}  Most of the information in this section
comes from \cite{Veech:eisenstein}. Let
  $$ \text{$Q_n$ denote the triangle with angles $\displaystyle
\frac{\pi}{n}, \frac{\pi}{n}, \frac{(n-2)\pi}{n}$,} $$
  (realized with the two equal sides having length $1$, and one of the
equal sides horizontal).  Then the surface corresponding to $Q_n$ can
easily be seen to be (isomorphic to) $X_n$.  The cylinder
decomposition in the vertical direction consists of $(n-1)/2$
cylinders $V_j$, and for $1 \le j \le (n-1)/2$, we have:
  \begin{gather}
\label{eq:heights:widths}
h_j = \text{ height $V_j$ } = 4 \sin \frac{\pi (2 j-1)}{n} \cos
\frac{\pi}{n} \\
w_j = \text{ width $V_j$ } = 2 \sin \frac{ \pi (2j-1)}{n} \sin 
\frac{\pi}{n}
\end{gather}
(the closed trajectories in the cylinder $V_j$ have length $h_j$).
Since for all $1 \le j \le (n-1)/2$, $h_j/w_j = 2 \cot \frac{\pi}{n}$,
the unipotent
\begin{displaymath}
u_n = \begin{pmatrix} 1 & 0 \\ 2 \cot \frac{\pi}{n} & 1 \end{pmatrix}
\end{displaymath}
belongs to the Veech group $\Gamma_n$ of $X_n$. Note that
\begin{equation}
\label{eq:move:perp}
u_n \begin{pmatrix} w_j \\ 0 \end{pmatrix} = \begin{pmatrix} w_j \\ h_j
\end{pmatrix}
  . \end{equation}
The unipotent $u_n$, together with the rotation by
$2 \pi/n$ generate $\Gamma_n$. It is shown in \cite{Veech:eisenstein} 
that
\begin{displaymath}
\Vol( \Gamma_n \backslash \half ) = \frac{n-2}{n} \pi
\end{displaymath}
where $\Vol$ denotes the Poincar\'e volume on the hyperbolic plane
$\half$.

The following lemma is from \cite{Gutkin:Judge}:
\begin{lem}
\label{lemma:GJ}
Suppose $\Gamma \subset \SL(2,\real)$ is a lattice, and suppose that
$\Gamma$ intersects non-trivially the stabilizer $N$ in
$\SL(2,\real)$ of $v$. {\rm(}The above condition is equivalent to the
discreteness of the orbit $\Gamma v$.{\rm)} Let $\gamma \in \Gamma$ be
either of the two generators of $\Gamma \cap N$. Let $B(T)$ denote the
ball in $\real^2$ of radius $T$ centered at the origin. Then, as $T
\to \infty$,
\begin{displaymath}
| \Gamma v \cap B(T) | \sim \Vol( \Gamma_n \backslash \half )^{-1}
\frac{ |\langle \gamma v^{\perp}, v \rangle |}{ \|v\|^3 \|v^\perp\|} 
T^2,
\end{displaymath}
where $v^\perp$ is any vector perpendicular to $v$.
\end{lem}
We also record the following trivial consequence (of the existence of
the asymptotics):

Suppose $v$ is as in Lemma~\ref{lemma:GJ}, and suppose $v'$ is a
scalar multiple of $v$. Then, as $T \to \infty$.
\begin{equation}
\label{eq:trivial}
|\Gamma v' \cap B(T)| \sim \frac{\|v\|^2}{\|v'\|^2} |\Gamma v \cap B(T)|
  . \end{equation}

We now apply Lemma~\ref{lemma:GJ} with $\Gamma = \Gamma_n$, $v =
\begin{pmatrix} 0 \\ h_1\end{pmatrix}$, $v^\perp = \begin{pmatrix} w_1
   \\ 0 \end{pmatrix}$, $\gamma = u_n$, and using
(\ref{eq:move:perp}) we get that the number of
cylinders in the ball of radius $T$ that are in the orbit of the
cylinder $V_1$ is asymptotic to
\begin{equation}
\label{eq:count:one:orbit}
\Vol( \Gamma_n \backslash \half )^{-1} \frac{ h_1^2}{h_1^3 w_1} T^2 =
\frac{n}{(n-2) \pi} \frac{1}{h_1 w_1} T^2
  . \end{equation}
We now use (\ref{eq:trivial}) to see that for $1 \le j \le (n-1)/2$,
the number of cylinders in
the ball of radius $T$ that are in the $\Gamma_n$ orbit of $V_j$ is
asymptotic to
\begin{equation}
\label{eq:count:j:one:orbit}
\frac{n}{(n-2) \pi} \frac{1}{h_j w_j} T^2.
\end{equation}
(in the above, we used the identity $\frac{1}{h_1 w_1}
\frac{h_1^2}{h_j^2} = \frac{1}{h_j w_j}$).
Since every cylinder
is in the orbit of some $V_j$, we get (after summing over $j$),
\begin{equation}
\label{eq:count:Veech:prelim}
N(X_n,T) \sim \frac{n}{(n-2) \pi} \left(\sum_{j=1}^{(n-1)/2} 
\frac{1}{h_j w_j}\right) T^2
\end{equation}
Using the identity \cite[Lemma 6.3]{Veech:eisenstein}
\begin{displaymath}
\sum_{j=1}^{(n-1)/2} \frac{1}{\sin^2 \frac{\pi (2j-1)}{n}} = \frac{n^2
     - 1}{6}
\end{displaymath}
and using the expressions (\ref{eq:heights:widths}), we get
\begin{displaymath}
N(X_n,T) \sim \frac{n}{(n-2) \pi} \frac{ n^2 -1 }{48} \frac{T^2}{\sin
   \frac{\pi}{n} \cos \frac{\pi}{n} }
  . \end{displaymath}
Since
\begin{displaymath}
\sin \frac{\pi}{n} \cos \frac{\pi}{n} = \frac{1}{2} \sin \frac{ 2
   \pi}{n} = \area Q_n = \frac{1}{2 n} \area X_n,
\end{displaymath}
we have
\begin{equation}
\label{eq:count:Veech:final}
N(X_n,T) \sim \frac{n^2(n^2-1)}{ 24 (n-2) \pi} \frac{T^2}{\area(X_n)}
= \frac{\pi}{\zeta(2)} \frac{n^2(n^2-1)}{144 (n-2)} \frac{T^2} 
{\area(X_n)}
  . \end{equation}
This is the formula in \cite{Veech:eisenstein}.

\subsection*{The Siegel-Veech formula applied to $X_n$.}
It is useful for the sequel to compare the result of
(\ref{eq:count:j:one:orbit}) with the result of the corresponding
Siegel-Veech formula. Let $D_n = \SL(2,\real) X_n$ denote the orbit of
$X_n$. This is a closed submanifold of the stratum, which is also
called a Teichm\"uller curve. For $1 \le j \le (n-1)/2$,
define $U_j \colon D_n \to \text{ subsets of } \real^2$
by  the formula
$U_j(g X_n) = g \Gamma_n \begin{pmatrix} 0 \\ h_j \end{pmatrix}$
(where $g \in \SL(2,\real)$).
Let $f_\epsilon$ denote the characteristic function of the
$\epsilon$-ball in $\real^2$ centered at the origin, and for
$M \in D_n$, define the Siegel-Veech transform
$\hat{f}_{j,\epsilon}(M) = \sum_{v \in U_j(M)} f_\epsilon(v)$.

\begin{lem}
\label{lemma:disjoint:cusps}
If $\epsilon$ is sufficiently small, then
$\hat{f}_{j,\epsilon} \colon D_n \to \real$
takes on only the values $0$ and $1$; we have
$\hat{f}_{j,\epsilon}(M) = 1$ if and only if $M$ has a cylinder
decomposition such that the $j^\text{th}$ cylinder from the left has
height at most~$\epsilon$. Given $M$, such a cylinder
decomposition is unique if it exists.
\end{lem}

\begin{proof}
  This is straightforward. (The uniqueness of the decomposition follows
from the fact, proved by Veech \cite{Veech:eisenstein}, that
$\half^2/\Gamma_n$ has only one cusp.)
  \end{proof}

Let $\nu$ denote the normalized
$\SL(2,\real)$ invariant measure on $D_n$. Then, we have the
Siegel-Veech formula:
\begin{equation}
\label{eq:Xn:V1:Siegel:Veech}
\int_{D_n} \hat{f}_{j,\epsilon} \,d \nu = c_j \int_{\real^2} f_\epsilon.
\end{equation}
We now apply Theorem~\ref{theorem:counting:EMS:EM} with $V(\cdot) =
U_j( \cdot)$, and $S = X_n$. The validity of assumption \pref{(C)}
can be deduced from the mixing property of the geodesic flow,
see \cite{Margulis:thesis} for a general proof in variable
negative curvature,  or \cite{Eskin:McMullen} for a simplified
exposition in the  constant curvature case. We obtain that
\begin{displaymath}
|U_j(X_n) \cap B(T)| \sim \pi c_j T^2,
\end{displaymath}
where $c_j$ is as in (\ref{eq:Xn:V1:Siegel:Veech}).
Comparing with (\ref{eq:count:j:one:orbit}), we see that
\begin{displaymath}
c_j = \frac{n}{(n-2) \pi^2} \frac{1}{h_j w_j}.
\end{displaymath}
Substituting into (\ref{eq:Xn:V1:Siegel:Veech}) we get
\begin{equation}
\label{eq:measure:small:V1}
\frac{1}{\pi \epsilon^2}
\int_{D_n} \hat{f}_{j,\epsilon} \,d \nu = \frac{n}{(n-2) \pi^2} 
\frac{1}{h_j
   w_j}.
\end{equation}

\begin{rmk}
It is possible to prove (\ref{eq:measure:small:V1})
directly, and thus to compute the asymptotics in 
(\ref{eq:count:one:orbit})
without using Lemma~\ref{lemma:GJ}. We chose this indirect derivation of
(\ref{eq:measure:small:V1}) to minimize the amount of computation.
\end{rmk}

\section*{The branched cover.}
We now return to our surface $S_n$, which is a branched cover of
$X_n$ (see Figure~\ref{fig:pentagons}). $X_n$ is a union of two
$n$-gons, and the two branch points $p$ and $p'$ are at the centers of 
the
$n$-gons.  We now wish to apply Theorem~\ref{theorem:largecircles} to
the point $(X_n,p,p') \in \cX^2$.

It is important to note that $X_n$ is hyperelliptic, and that our
two branch points are interchanged by the hyperelliptic
involution. Since the hyperelliptic involution commutes with the
$\SL(2,\real)$ action, this it true for any point in the orbit of 
$(X_n,p,p')$.
Thus, the $\SL(2,\real)$ orbit of $(X_n,p,p')$
is {\em not} dense in the space $\cX^2$, and indeed we have in
Theorem~\ref{theorem:largecircles} a proper $W \subset (\real^2)^2$ of
real dimension $2$.
Let $L$ denote the subspace $\{\,(v, -v),  \mid v \in \real^2\,\}$.
The above argument shows that $W \subseteq L$. But since we know that
$S_n$ is not Veech, $\dim W > 0$. Hence $\dim W = 2$ and $W = L$.
Let $\cO = (\SL(2,R) \ltimes \Phi^2_W) (X_n,p,p')$. Then $\cO \subset
\cX^2$
consists of points of the form $(M, q,q')$ where $M \in D_n$, $q \in
M$, $q' \in M$ and $q$ and $q'$ are interchanged by the hyperelliptic
involution of $M$. By Theorem~\ref{theorem:largecircles},
\begin{equation}
\label{eq:marked:points}
\lim_{t\to \infty} \nu_{t,(X_n,p,p')} = \mu,
\end{equation}
where $\mu$ is Lebesgue measure on $\cO$.

Now let $\tilde{\cO}$ denote the orbit closure $\overline{\SL(2,\real)
   S_n}$. Since $S_n$ is a double cover of $X_n$, branched over $p$ and
$p'$, for any $g \in \SL(2,\real)$, $g S_n$ is a double cover of $g
X_n$ branched over $gp$ and $gp'$, and $(g X_n, gp, gp')\in \cO$.
Thus, in particular, every surface in $\tilde{\cO}$ is a double cover
of a surface in $D_n$.
Thus we have a natural map $\tilde{\pi} \colon \tilde{\cO} \to \cO$
that maps each surface $S \in \tilde{\cO}$ to the surface in $D_n$ of
which it is a double cover, and notes the locations of the branch
points.  Now in view of (\ref{eq:marked:points}) and
Lemma~\ref{FinCovClassify},
\begin{displaymath}
\lim_{t\to \infty} \nu_{t,S_n} = \tilde{\mu},
\end{displaymath}
where $\tilde{\mu}$ is normalized Lebesgue measure on $\tilde{\cO}$.
Hence, by Theorem~\ref{theorem:counting:EMS:EM}, we have a quadratic
asymptotic formula
\begin{displaymath}
N(S_n,T) = |V_{cyl}(S_n) \cap B(T)| \sim \pi c T^2,
\end{displaymath}
with the constant $c$ given by
\begin{equation}
\label{eq:c:tilde:O}
c = \frac{1}{\pi \epsilon^2} \int_{\tilde{\cO}} \hat{f}_\epsilon \,
d\tilde{\mu}
\end{equation}
where as above, $f_\epsilon \colon \real^2 \to \real$ is the
characteristic function of the ball of radius $\epsilon$ centered at
the origin, and  $\hat{f}_\epsilon(S) = \sum_{ v \in V_{cyl}(S)}
f_\epsilon(v)$.

Let $v$ be some periodic direction for $S_n$, hence for $X_n$. We may
use an element $\gamma$ of the Veech group $\Gamma_n$ of $X_n$
to map $v$ to the vertical direction. Note that $\gamma S_n$ is a
double cover of $\gamma X_n = X_n$. In the vertical
direction, $X_n$ has the cylinder decomposition $V_1, \dots,
V_{(n-1)/2}$ described above.
\begin{lem}
\label{lemma:cyl:decomp}
For any $\gamma \in \Gamma_n$, the branch points of $\gamma S_n$ will
project to two points in the same cylinder, say $V_k$.
The cylinder decomposition of $\gamma S_n$ in the vertical
direction is the following:
\begin{enumerate} \renewcommand{\theenumi}{\alph{enumi}}
\item \label{lemma:cyl:decomp-nobranch}
  For each $j \ne k$, there are two cylinders on
   $\gamma S_n$ of the same length as $V_j$ {\rm(}one on each
``sheet''{\rm)}.  \item \label{lemma:cyl:decomp-branch}
  On $\gamma S_n$ there are two cylinders of
the same length as $V_k$ and two cylinders of twice the length of
$V_k$.
\end{enumerate}
\end{lem}

\begin{proof} The fact that both branch points project to the same
cylinder of $X_n$ follows from the fact that each cylinder of $X_n$ is
preserved by the hyperelliptic involution $\sigma$ of $X_n$ (since
different cylinders have different lengths) and the fact that the
branch points are interchanged by $\sigma$.   From
Figure~\ref{fig:pentagons}, the cover $S_n$ is determined by two slits
(drawn as the thick lines in the figure), which are interchanged by
$\sigma$. Since $\sigma$  commutes with the $\SL(2,\real)$ action, the
cover $\gamma S_n$ of $X_n$ is also determined by two slits, which are
interchanged by $\sigma$.  For each cylinder $V_j$ of $X_n$, let
$\lambda_j$ denote the closed trajectory in the center of $V_j$. Note
that for any $j$, $\lambda_j$ is mapped to itself under $\sigma$.
Also, since $\sigma$ exchanges the slits, we see that $\lambda_j$
intersects each slit the same number of times. Thus, $\lambda_j$
breaks up into two closed paths of the same length when lifted from
$X_n$ to $\gamma S_n$. This proves \pref{lemma:cyl:decomp-nobranch}
and the first assertion of \pref{lemma:cyl:decomp-branch}.  It is easy
to see that the closed vertical trajectories on $V_k$ between the
boundary of $V_k$ and one of the branch points double in length when
lifted from $X_n$ to $\gamma S_n$. This proves the second assertion of
\pref{lemma:cyl:decomp-branch}.
  \end{proof}

\begin{cor}
\label{cor:descent} The function $\hat{f}_\epsilon \colon \tilde{\cO}
\to \real$ is constant on the fibers of $\tilde{\pi}$ {\rm(}a.e.{\rm)},
and thus descends to a function $\bar{f}_\epsilon \colon \cO \to
\real$. The latter function, for $\epsilon$ sufficiently small,  is
given {\rm(}a.e.{\rm)} by the formula
  \begin{equation} \label{cor:descent:eqn}
\bar{f}_\epsilon(M,q,q') = 2 \hat{f}_{k,\epsilon}(M) +
2 \hat{f}_{k,\epsilon/2}(M) + 2 \sum_{j \ne k}
\hat{f}_{j,\epsilon}(M),
\end{equation}
where for $1 \le j \le (n-1)/2$, $\hat{f}_{j,\epsilon} \colon D_n \to
\real$ is as in Lemma~\ref{lemma:disjoint:cusps}, and
$k$ is such that $q$ {\rm(}and $q'${\rm)} belong to the $k^\text{th}$
cylinder from the left in the unique cylinder decomposition that
contains a cylinder of height at most $\epsilon$.
\end{cor}

\begin{proofsketch}
  Choose a fundamental domain for~$\Gamma_n$ in the upper half plane.
Since $X_n$ has a vertical cylinder decomposition, we may assume that
the cusp of the fundamental domain approaches~$\infty$, rather than
approaching a point on the real axis. This means that as $h$ goes to
infinity in the fundamental domain, the unique short cylinder
decomposition of $h X_n$ is the image under~$h$ of the vertical
cylinder decomposition of~$X_n$.

We first prove that \pref{cor:descent:eqn} is
correct for all $M$ in the $\SL(2,\real)$ orbit of~$S_n$. To do
this, let $g \in \SL(2,\real)$, and write $g = h \gamma$,
where $\gamma \in \Gamma_n$, and $h$ is in the fundamental
domain.
  Note that $\hat{f}_\epsilon(M)$ is zero unless $M$ has a short
cylinder. Thus, if $g = h \gamma$, and $h$ is in a compact part
of the fundamental domain, then (in view of
Lemma~\ref{lemma:cyl:decomp}), we have $\hat{f}_\epsilon(h
\gamma S_n) = 0$. Therefore, we may assume that $h$ is in the
cusp, and hence the unique short cylinder decomposition of $g S_n = h
\gamma S_n$ is the image, under the linear action of~$h$, of
the vertical cylinder decomposition of $\gamma S_n$. Then it is clear
from Lemma~\ref{lemma:cyl:decomp} that \pref{cor:descent:eqn} holds
for $M = g S_n$.

To complete the proof, note that both sides of
\pref{cor:descent:eqn} are continuous off a closed set of
measure~$0$, namely, the set where a branch point projects to
an edge of a cylinder (in a cylinder decomposition containing a
cylinder of height at most~$\epsilon$). Then use the fact that the
$\SL(2,\real)$ orbit of~$S_n$ is dense.
  \end{proofsketch}

\begin{rmk}
The analogue of the first assertion
of Corollary~\ref{cor:descent} fails in the context of
\cite{EMS}, in part since there we are dealing with covers of high
degree. This is responsible for most of the combinatorial complexity
of the argument in \cite{EMS}.
\end{rmk}
\medskip

In view of Corollary~\ref{cor:descent}, (\ref{eq:c:tilde:O}) becomes
\begin{align*}
c & = \frac{1}{\pi \epsilon^2} \int_{\cO} \bar{f}_\epsilon \, d \mu \\
   & = \frac{1}{\pi \epsilon^2} \int_{D_n} \int_M
   \bar{f}_\epsilon(M,q,\sigma(q)) \, d\lambda_M(q) \, d\nu(M),
\end{align*}
where $\lambda_M$ is Lebesgue measure on the translation surface $M$, 
and
$\sigma$ denotes the hyperelliptic involution.
Performing the integral over $M$, we get
\begin{displaymath}
   c = \sum_{k=1}^{(n-1)/2} p_k \left( \frac{1}{\pi \epsilon^2}
     \int_{D_n} \left( 2 \hat{f}_{k,\epsilon}(M) +
2 \hat{f}_{k,\epsilon/2}(M) + 2 \sum_{j \ne k}
\hat{f}_{j,\epsilon}(M)\right) \, d\nu(M) \right),
\end{displaymath}
where $p_k = h_k w_k/A$ and $A = \area(X_n)$ (so that $p_k$ denotes
the relative area of the $k^\text{th}$ cylinder from the left in any
cylinder decomposition).

Now, using (\ref{eq:measure:small:V1}), we get
%
%
%
\begin{align*}
N(S_n,T)   \sim c \, T^2 & = \frac{n}{(n-2) \pi} \sum_{k=1}^{(n-1)/2} 
p_k
\left(\frac{2}{h_k w_k} + \frac{2}{4 h_k w_k} + \sum_{j\ne k}
   \frac{2}{h_j w_j} \right) T^2 \\
& = \frac{n}{(n-2) \pi} \sum_{k=1}^{(n-1)/2} p_k
\left(\frac{2}{4 h_k w_k} + \sum_{j=1}^{(n-1)/2}
   \frac{2}{h_j w_j} \right) T^2 . \\
\end{align*}
Since the second term in the parenthesis is independent of $k$, and
$\sum p_k =1$, this can be rewritten as
\begin{equation}
N(S_n,T) \sim \left( \frac{n}{(n-2) \pi} \sum_{j=1}^{(n-1)/2}
   \frac{2}{h_j w_j} + \frac{n}{(n-2) \pi} \sum_{k=1}^{(n-1)/2} \frac{
     2 p_k}{ 4 h_k w_k} \right) T^2
  . \end{equation}
The first term in the parenthesis is in view of
(\ref{eq:count:Veech:prelim}) twice the limit of $N(X_n,T)/T^2$. The
second term in the parenthesis is, since $p_k = h_k w_k/(\area X_n)$,
equal to
\begin{displaymath}
\frac{n}{(n-2) \pi} \frac{(n-1)}{4 \area(X_n)}
  . \end{displaymath}
In view of (\ref{eq:count:Veech:final}), we get
\begin{displaymath}
N(S_n,T) \sim \left(\frac{n^2(n^2-1)}{ 12 (n-2) \pi}
   \frac{1}{\area(X_n)} + \frac{n}{(n-2) \pi} \frac{(n-1)}{4
     \area(X_n)} \right) T^2
  . \end{displaymath}
Simplifying, we get
\begin{displaymath}
N(S_n,T) \sim \frac{n(n-1)(n^2+n+3)}{12(n-2)\pi} \frac{T^2}{\area(X_n)}
  . \end{displaymath}
Alternatively,
\begin{displaymath}
N(S_n,T) \sim \frac{n(n-1)(n^2+n+3)}{6(n-2)\pi} \frac{T^2}{\area(S_n)}
= \frac{\pi}{\zeta(2)} \frac{n(n-1)(n^2+n+3)}{36(n-2)}
\frac{T^2}{\area(S_n)}
  . \end{displaymath}
   To pass from~$S_n$ to~$P_n$, note that $N(P_n,T) = N(S_n,T)$,
and since $S_n$ consists of $4 n$ triangles, $\area(S_n) = 4n
\area(P_n)$.

\acks This research was partially supported by grants from the National
Science Foundation (DMS-0244542 and DMS-0100438), the Packard
Foundation, and an EPSRC Advanced Research Fellowship.
  The authors are grateful to the Swiss Federal Technical Institute (ETH)
of Zurich, the University of Chicago, the Newton Institute (Cambridge,
England), and the American Institute of Mathematics for their
hospitality, and to Kariane Calta, Fran\c cois Ledrappier, Howard
Masur, Martin Schmoll, and anonymous referees for helpful comments on a
preliminary version of the manuscript.

\end{document}